\documentclass{amsart}
\usepackage{amsfonts,amssymb,amscd,amsmath,enumerate,verbatim,calc}

\newcommand{\CM}{Cohen-Macaulay}

\newcommand{\n}{\mathfrak{n} }
\newcommand{\m}{\mathfrak{m} }

\newcommand{\p}{\mathfrak{p} }

\newcommand{\spec}{\operatorname{Spec}}
\newcommand{\mspec}{\operatorname{m-Spec}}

\newcommand{\rank}{\operatorname{rank}}

\newcommand{\Ass}{\operatorname{Ass}}

\newcommand{\ann}{\operatorname{ann}}
\newcommand{\grade}{\operatorname{grade}}
\newcommand{\depth}{\operatorname{depth}}

\newcommand{\Hom}{\operatorname{Hom}}

\theoremstyle{plain}

\newtheorem{theorem}{Theorem}[section]
\newtheorem{corollary}[theorem]{Corollary}
\newtheorem{lemma}[theorem]{Lemma}
\newtheorem{proposition}[theorem]{Proposition}

\theoremstyle{definition}
\newtheorem{definition}[theorem]{Definition}
\newtheorem{s}[theorem]{}
\newtheorem{remarks}[theorem]{Remarks}
\newtheorem{remark}[theorem]{Remark}
\newtheorem{example}[theorem]{Example}
\newtheorem{examples}[theorem]{Examples}

\begin{document}
\title[Ratliff-Rush filtration]{Ratliff-Rush Filtrations Associated With Ideals and Modules over a Noetherian Ring}
\author{Tony~J.~Puthenpurakal}
\address{Department of Mathematics, IIT Bombay, Powai, Mumbai 400 076, India}
\email{tputhen@math.iitb.ac.in}
\author{Fahed~Zulfeqarr}
\address{Department of Mathematics, IIT Bombay, Powai, Mumbai 400 076, India}
\email{fahed@math.iitb.ac.in}
\thanks{The first author was partly supported by IIT Bombay seed grant 03ir053. The second author thanks  CSIR  for financial support.}

\begin{abstract}
Let $R$ be a commutative Noetherian ring, $M$ a finitely generated $R$-module and $I$ a proper ideal of $R$. In this paper we introduce and  analyze some properties of $r(I, M)=\bigcup_{k\geqslant 1} (I^{k+1}M: I^kM)$, {\it the Ratliff-Rush ideal associated with $I$ and $M$}. When $M= R$ (or more generally when $M$ is projective) then
$r(I, M)= \widetilde{I}$, the usual Ratliff-Rush ideal associated with $I$.
 If $I$ is a regular ideal and $\ann M=0$ we show that $\{ r(I^n,M) \}_{n\geqslant 0}$ is a stable $I$-filtration. 
If $M_{\p}$ is free for all ${\p}\in \spec R\setminus \mspec R,$  then under mild condition on $R$ we show that
 for a regular ideal $I$,    $\ell( r(I,M)/{\widetilde I})$ is finite. Further $r(I,M)=\widetilde I\;$ if  $A^*(I)\cap \mspec R =\emptyset $  (here  $A^*(I)$ is the stable value of the sequence $\Ass (R/{I^n})$). Our generalization also
 helps to better understand the usual Ratliff-Rush filtration. When $I$ is a regular $\m$-primary ideal our techniques yield an easily computable bound for $k$ such that $\widetilde{I^n} = (I^{n+k} \colon I^k)$ for all $n \geqslant 1$. For any
 ideal $I$ we show that $\widetilde{I^nM}=I^nM+H^0_I(M)\quad\mbox{for all}\;n\gg 0.$ This  yields that $\widetilde {\mathcal R}(I,M)=\bigoplus_{n\geqslant 0} \widetilde {I^nM}$ is Noetherian if and only if $\depth M>0$. Surprisingly if $\dim M=1$ then $\widetilde G_I(M)=\bigoplus_{n\geqslant 0} \widetilde{I^nM}/{\widetilde{I^{n+1}M}}$ is always a Noetherian and a Cohen-Macaulay $G_I(R)$-module. Application to Hilbert coefficients is also discussed.
\end{abstract}
\date{\today}
\maketitle
\section*{Introduction}
\noindent Let $R$ be a commutative Noetherian ring and $I$ an ideal of $R$. The Ratliff-Rush ideal $\widetilde{I} = \bigcup_{k\ge 1} (I^{k+1}\colon I^k)$ is a useful notion. When $R$ is local and $I$ is $\m$-primary, the Ratliff-Rush filtration $\{ \widetilde{I^n} \}_{n\geq 1}$ has many applications in the theory of Hilbert functions, for instance
see \cite{RV2}.
 In this paper we generalize this notion.\\
\indent  Let  $M$ be a finitely generated $R$-module. Set
$$r(I, M)=\bigcup_{k\ge 1} (I^{k+1}M: I^kM).$$
We call  $r(I,M)$ the {\it Ratliff-Rush ideal associated with $I$ and $M$}. Notice that $r(I,M)=\widetilde I$ if $M=R$.
Our generalization also gives us a better understanding of the usual Ratliff-Rush  filtration. For instance when  $R$ is a \CM \ local ring  and $I$ is an $\m$-primary ideal  it is  useful  find a upper bound on $k$ such that
$\widetilde{I^n} = (I^{n+k} \colon I^k)$ for all $n \geqslant 1$, see \cite{RV}.
Our techniques enables us to find an easily computable upper bound on $k$
(see  \ref{21} and  \ref{obs}).

\indent We analyze many of its properties. Perhaps the first non-trivial property is  that the function $\;I\longmapsto r(I,M)\;$ is an involution on the set of ideals of $R$,  i.e., 
 $$r(r(I, M), M) =  r(I, M).$$
This is done in Theorem \ref{closedness2}. 

Next we relate this notion to integral closure. In Theorem \ref{10} we show that if $I$ is a regular ideal then there exists a
rank $1$ module $M$ such that $r(I, M) = \overline{I}$. A typical example of a rank one module is a regular ideal. In Proposition \ref{rank2} we show that there exists a regular ideal $J$ such that  $r(I,J) = \overline{I}$. Furthermore
we  prove (in Proposition \ref{rank2}) that the
set
\[
\mathcal{C}(I) = \{ J \mid J \ \text{a regular ideal with }  \ r(I, J) = \overline{I}  \}.
\]
has a unique maximal element.\\
\indent Next we analyze the filtration $ \mathcal{F}_{M}^{I} = \{r(I^n, M)  \}_{n\geq 1}$. We first prove that this is a filtration of ideals and an $I$-filtration (see Theorem \ref{1}).  Thus $\;{\mathcal R}({\mathcal F}^I_M) = \bigoplus_{n \geq 0}  r(I^n, M) $ is a $\mathcal{R}(I) = \bigoplus_{n \geq 0}I^n$-algebra. 
 Let $\widetilde{IM}=\cup_{k\ge 1} (I^{k+1}M:_M I^k)$, the Ratliff-Rush module of $M$ associated with $I$.
We show that  ${\widetilde {\mathcal R}}(I,M) = \bigoplus_{n \geq 0} \widetilde{I^nM}$ is a graded  $\;{\mathcal R}({\mathcal F}^I_M)$-module (see  Proposition \ref{6}).
In Theorem \ref{stablity1} we prove that if $\grade(I, M) > 0\;$ and $\;\ann M = 0\;$ then $\; \mathcal{F}_{M}^{I}$ is a stable $I$-filtration.

If $M$ is a projective $R$-module then $r(I,M) = \widetilde{I}$ for all ideals $I$ (see Proposition \ref{premproj}).
As $R$ is Noetherian, projective (finitely generated) modules are precisely locally free (finitely generated) modules. 
In some sense the next case is to consider $R$-modules $M$ such that
\begin{equation*}
 M_{\p} \ \text{is free} \  \text{for all} \ \p \in \spec(R) \setminus \mspec(R) \tag{*}
\end{equation*}
(here $\mspec R$ denotes the set of maximal ideals of $R$). For instance if $M \subset F$ where $F$ is a free $R$-module
and if $\ell(F/M)$ is finite then $M$ satisfies (*). Here $\;\ell (-)$ denotes the length. On $R$ we impose a mild condition $\Ass R \bigcap \mspec R =\emptyset $ (every domain that is not a field satisfies this condition).
 We show that if $I$ is a regular ideal then $\ell\left( r(I,M)/{\widetilde I}\right)$ is finite and the function $n\longmapsto \ell( r(I^n,M)/{\widetilde I^n})$ is a polynomial function. We show that $r(I,M)=\widetilde I\;$ if  $\;A^*(I)\cap \mspec R =\emptyset $. Here $A^*(I)$ is  the stable value of $\Ass(R/{I^n})$.\\
\indent In the case when $\ann M \neq 0$ or $\grade(I,M) = 0$  we go modulo the ideal $(H^{0}_{I}(M) \colon M)$. In Proposition \ref{recursion}  we show that if $M \neq  H^{0}_{I}(M)$ then 
\[
r(I^{n+1}, M) = I\cdot r(I^{n}, M) + (H^{0}_{I}(M) \colon M) \quad \text{for all} \ n \gg 0.
\]
An easy consequence of our techniques (see Proposition \ref{impcor1}) is the following result  in the case of the usual Ratliff-Rush filtration of a module $M$ with respect to $I$
\[
\widetilde{I^nM} = I^nM + H^{0}_{I}(M)  \quad \text{for all} \ n \gg 0.
\]

\indent In the final section we show that if $\dim M=1$, $\widetilde G_I(M)=\bigoplus_{n\geqslant 0} \widetilde{I^nM}/{\widetilde{I^{n+1}M}}$ is a Noetherian and a Cohen-Macaulay $G_I(R)$-module. This is surprising since if $\depth M=0$ (and $\dim M=1$) then $\widetilde {\mathcal R}(I,M)=\bigoplus_{n\geqslant 0} \widetilde {I^nM}$ is not a Noetherian ${\mathcal R}(I)$-module. We also give an application of our result to Hilbert coefficients.

\indent Here is an overview of the contents of the paper. In section one we study few basic properties of the
ideal $r(I,M)$. In section two we study the filtration $\{r(I^n, M)\}_{n\geq 1}$ and explore the relation between $\{r(I^n, M)  \}_{n\geq 1}$ and $\{\widetilde{I^nM}\}_{n\ge 1}$. In section three we prove that the operation $I \mapsto r(I,M)$ is an involution. In section four we relate it to integral closure. In section five we prove that it is a stable $I$-filtration when $I$ is  regular and $\ann M = 0$. In  section six we study the case when $M_{\p}$ is free for all $\p\in \spec R\setminus \mspec R$.   In section seven we study the general case when $\ann M \neq 0$ or $\grade(I,M) = 0$. For the next sections we assume that $(R,\m)$ is local. In section eight we study its relation with superficial elements. This is then used to give a bound on $k$ such that $\widetilde{I^n} = (I^{n+k} \colon I^k)$ for all $n \geq 1$. In section nine we study ideals having a principal reduction and use it to compute  $r(I,M)$ in some examples. In final section we show that if $\dim M=1$ then  $\widetilde G_I(M)=\bigoplus_{n\geqslant 0} \widetilde{I^nM}/{\widetilde{I^{n+1}M}}$ is always a Cohen-Macaulay $G_I(R)$-module. Finally we give an application of one of our  results to Hilbert coefficients of a 1-dimensional module.

\section{Preliminaries}
\noindent In this paper unless otherwise stated all rings considered are commutative Noetherian and all modules are assumed to be finitely generated.\\ 
\indent Let $R$ be a ring and $I$ an ideal of $R$. Let $M$ be an $R$-module. Consider the following ascending chain of ideals in $R$
$$ I\subseteq (IM : M)\subseteq (I^2M : IM)\subseteq ...\subseteq (I^{k+1}M : I^kM)\subseteq ...$$
Since $R$ is Noetherian, this chain stabilizes. We denote its stable value by $r(I,M)$.  We call  $r(I,M)$ the {\it Ratliff-Rush ideal associated with $I$ and $M$}.\\
\indent In this section we prove some basic properties of ideal $\;r(I,M)$, in particular we show that 
$$r(I,M\oplus N)=r(I,M) \cap r(I,N)\quad \mbox{and}\quad r(I,M\otimes_R N) \supseteq r(I,M) + r(I,N).$$
We also investigate the case when $r(I,M)=R.$

\begin{remarks}\label{firstremarks}$~~~$
\begin{itemize}
\item [$($\rm a$)$] $r(I,M) =(I^{k+1}M : I^kM)\quad \text{for all}\;k\gg 0.$
\item [$($\rm b$)$] When $M=R$, 
$$r(I, R)= \bigcup_{k\ge 1} (I^{k+1}: I^k)=\widetilde{I},\;\mbox{the Ratliff-Rush closure of}\; I.\quad\mbox{(see \cite{RR})}$$
\item [$($\rm c$)$] For $n\ge 1$, we have 
$$r(I^n,M)=\bigcup_{k\ge 1} (I^{n+k}M: I^kM).$$
\item [$($\rm d$)$] One can easily check that
$$I\subseteq \widetilde{I} \subseteq r(I, M).$$
\item [$($\rm e$)$]
Let $R=\bigoplus_{n\ge 0}R_n$ be a graded ring and $I$, a homogeneous ideal of $R$. Let $M$ be a graded $R$-module. Then $r(I,M)$ is a homogeneous ideal.
\end{itemize}
\end{remarks}

We give an example which shows that there exists a module $M$ such that
$$I\subsetneq{\widetilde I} \subsetneq r(I,M).$$ 

\begin{example}\label{9} 
Let $R= k[ t^4,t^{11},t^{17},t^{18}]$,\quad $\m=\left\langle t^4,t^{11},t^{17},t^{18}\right\rangle R$,\quad $I=\left\langle t^4,t^{11}\right\rangle R$\quad and\; $M=\left\langle t^4,t^{11},t^{17}\right\rangle R$.
One can prove (by induction on $n$) that
\begin{eqnarray*}
I^n &=& \left\langle t^{4n+7i}\; : i=0,1,2,...,n\right\rangle R\quad\mbox{for all} \; n\ge 1. 
\end{eqnarray*}
Note that $t^{18}\cdot I^2 \subseteq I^3\;$ so $t^{18}\in  {\widetilde I}$. One can verify that $\;t^{17} \notin (I^{n+1}: I^n)$ for all $n\ge 1$. Thus we get $\;{\widetilde I}=\;\left\langle t^4,t^{11},t^{18}\right\rangle R$.  Notice that $\;t^{17} \in (IM:M)\;$. Now since $\widetilde I\subseteq r(I,M)$, we get $r(I,M)=\left\langle t^4,t^{11},t^{17},t^{18}\right\rangle R$.
Therefore $I\subsetneq {\widetilde I}\subsetneq r(I,M) $ and  $r(I, M)= \m$. 
\end{example}

We next give an example from \cite[1.4]{RS} of a Ratliff-Rush closed ideal $I$  of ring $R$  and a module $M$ such that $I \ne r(I,M)$.

\begin{example}\label{equality}
Let $R=k[x,y]$ be a polynomial ring in variables $x$ and $y$. Let 
\begin{eqnarray*}
I &=& \; \left\langle y^{22},x^4y^{18},x^7y^{15},x^8y^{14},x^{11}y^{11},x^{14}y^8,x^{15}y^7,x^{18}y^4,x^{22}\right\rangle R,\\
\mbox{and}\;\;\;M &=& \left\langle x^2,y^3\right\rangle R.
\end{eqnarray*}
By \cite[1.4]{RS}, $I$ is Ratliff-Rush closed, that is $I=\widetilde I$ (also see Example \ref{exampleone}). Using {\it Singular} \cite{Sin} one can  verify that
\begin{eqnarray*}
(I^4M:I^3M) &=& I\;+\left\langle x^2y^{21},x^6y^{17},x^{13}y^{10},x^{20}y^3\right\rangle R.
\end{eqnarray*} 
So, $I\subsetneq (I^4M:I^3M)\subseteq r(I,M).$
Hence $I=\widetilde I \subsetneq r(I,M)$.
\end{example}

We collect some properties of the ideal $r(I,M)$ in the following proposition.

\begin{proposition}\label{3}
 Let $M$ and $N$ be any two $R$-modules. Then
\begin{itemize}
\item [$($\rm a$)$] If $f:M\longrightarrow N$ is a surjective $R$-homomorphism then
$r(I,M) \subseteq r(I,N)$.
\item [$($\rm b$)$] 
$r(I,M\oplus N)=r(I,M) \cap r(I,N).$
\item [$($\rm c$)$]
$r(I,M\otimes_R N) \supseteq r(I,M) + r(I,N).$
\item [$($\rm d$)$] If $T$ is a Noetherian ring which is a flat $R$-algebra then
$$r(I,M) \otimes_R T=r_T(IT,M\otimes_R T).$$
\item [$($\rm e$)$] If $S$ is a multiplicative closed subset of $R$ then\;
$$r(I,M)\otimes_R R_S= r_{R_S}(IR_S,M_S).$$
\item [$($\rm f$)$] For a regular ideal $I$,\;\
$r({\widetilde I},M)=r(I,M).$
\item [$($\rm g$)$]  For each $n\ge 1$, we have \;\;$r(I^n,M) = r(I^n,I^sM)\quad \mbox{for all}\;s\ge 1.$
\end{itemize}
\end{proposition}
\begin{proof} {\bf (a)} 
 Let $x\in r(I,M)$. Then, for some $k\ge 0$, we have $x \in  (I^{k+1}M:I^kM).$ Thus
\begin{eqnarray*}
 xI^kM &\subseteq & I^{k+1}M,\\
\mbox{so}\;\;\;\;xI^kf(M) &\subseteq & I^{k+1}f(M).
\end{eqnarray*}
Since $f$ is surjective, we get $x \in (I^{k+1}N:I^kN)\subseteq r(I,N).$\\
\noindent {\bf (b)} Let $x\in r(I,M\oplus N)$. Then $xI^k(M\oplus N)\subseteq I^{k+1}(M\oplus N)$ for some $k\ge 0$. Therefore 
\begin{eqnarray*}
xI^kM\subseteq I^{k+1}M\quad &\mbox{and}& \quad xI^kN\subseteq I^{k+1}N.\\
\mbox{So}\;\;\;x\in r(I,M)\cap r(I,N).
\end{eqnarray*}
The reverse inclusion is obvious.\\
\noindent {\bf (c)} Let $x\in r(I,M)$. We have $xI^kM\subseteq I^{k+1}M$ for some $k\ge 0$. Let $\sum_i \alpha_i (m_i\otimes n_i) \in I^k(M\otimes N)$, where $\alpha_i \in I^k$. We have
$$x\left( \sum_i\alpha_i (m_i\otimes n_i)\right)  = \sum_i (x\alpha_i m_i)\otimes n_i$$
Since $x\alpha_i m_i \in I^{k+1}M,\; x\alpha_i m_i=\sum_j\beta^{(i)}_j u^{(i)}_j$, where $\beta^{(i)}_j\in I^{k+1}$ and ${u}^{(i)}_j\in M$. Thus
\begin{eqnarray*}
x\left( \sum_i\alpha_i (m_i\otimes n_i)\right)  &=& \sum_i\left( \sum_j\beta^{(i)}_j{u}^{(i)}_j\right) \otimes n_i\\
&=& \sum_{i,j}(\beta^{(i)}_j{u}^{(i)}_j)\otimes n_i\\
&=& \sum_{i,j}\beta^{(i)}_j({u}^{(i)}_j\otimes n_i) \in I^{k+1}(M\otimes N).
\end{eqnarray*}
In the same way one can show that $r(I,N)\subseteq r(I,M\otimes_R N)$.\\ 
\noindent {\bf (d)} Using \cite[18.1]{Naga}, we have 
\begin{eqnarray*}
(I^{n+1}M:I^nM)\otimes_R T &=& (I^{n+1}M\otimes_R T:_T I^nM\otimes_R T)\quad \mbox{for all}\;n\ge 1.
\end{eqnarray*}
Therefore we get
\begin{eqnarray*}
 r(I,M) \otimes_R T &=& (I^{k+1}M:I^kM)\otimes_R T\;\;\mbox{for all}\;k\gg 0,\\
&=& (I^{k+1}M\otimes_R T:_T I^kM\otimes_R T)\;\;\mbox{for all}\;k\gg 0,\\
&=& r_T(IT,M\otimes_R T).
\end{eqnarray*}
\noindent {\bf (e)} This follows from part (d), since $R_S$ is a flat $R$-algebra.\\
\noindent {\bf (f)} For any regular ideal $I$ of $R$, we have  
$I^k= {\widetilde I}^k$ for all $k\gg 0$ (see \cite[2.1]{RR}). Therefore 
\begin{eqnarray*}
 r({\widetilde I},M) &=& ({\widetilde I}^{k+1}M:{\widetilde I}^kM) \quad\mbox{for all}\;k\gg 0,\\
 &=& (I^{k+1}M:I^kM) \quad\mbox{for all}\;k\gg 0,\\
&=& r(I,M).
\end{eqnarray*}
\noindent {\bf (g)} Notice that 
\[r(I^n,I^sM)=\bigcup_{k\geqslant 1}(I^{n+k}I^sM:I^kI^sM)=r(I^n,M).\]
\end{proof}

\begin{remark}\label{prelifreemodule} 
From  \ref{3}(b) and \ref{firstremarks}(d) it follows that if a module $M$ has a free summand
then $r(I,M)=\widetilde I$.
\end{remark}
\begin{corollary}\label{premproj}
Let $M$ be  a  projective $R$-module and let $I$ be an ideal of $R$. Then $r(I,M)=\widetilde I$.
\end{corollary}
\begin{proof}
Set $J=r(I,M)$. Since $M$ is projective, $M_{\p}$ is free for all $\p \in\spec R$, by \cite[4.11 (b)]{Eis}. Note that $\widetilde I \subseteq J$. Now if $\p \in\spec R$ then 
\begin{eqnarray*}
 J_{\p} &=& r_{R_{\p}}(I_{\p},M_{\p}), \;\;\mbox{(by \ref{3}e )}\\
\mbox{Also}\quad ( \widetilde I)_{\p} &=& \widetilde{I_{\p}}. 
\end{eqnarray*}
Thus $\;J/{\widetilde I}\;$ is locally zero and hence zero.
\end{proof}

The next proposition enables us to determine when is $r(I,M)=R$. Let 
$$H^0_I(M)=\{m\in M : I^nm=0\;\; \mbox{for some} \;n\ge 1\}$$
be the $I$-torsion submodule of $M$.

\begin{proposition}\label{basicreq}
Let  $I$ be a proper ideal of ring $R$. The following conditions are equivalent
\begin{itemize}
\item [$($\rm a$)$] $r(I,M)=R$.
\item [$($\rm b$)$] there exists $n\ge 0$ such that $I^nM=I^{n+1}M$.
\end{itemize}
If $H^0_I(M)=M$ then $r(I,M)=R$. Furthermore if  $R$ is local then the converse is also true.
\end{proposition}
\begin{proof}
 (a) $\Rightarrow $ (b): If $r(I,M)=R$ then $1\in (I^{n+1}M: I^nM)$ for some $n$. Therefore we get
$I^{n+1}M=I^nM.$

(b) $\Rightarrow$ (a): If $I^nM=I^{n+1}M$ for some $n\ge 0$ then $1\in (I^{n+1}M:I^nM)\subseteq r(I,M)$. So $r(I,M)=R$.\\
\indent  Suppose $H^0_I(M)=M$. Since $M$ is a finitely generated $R$-module, there exists an integer $n\ge 1$ such that $I^nM=0$. Hence $r(I,M)=R$.\\
\indent Let $R$ be local and $I\subseteq \m$. If $r(I,M)=R$ then by (b), $I^nM=I^{n+1}M$. So by
Nakayama's lemma,  $I^nM=0$. Thus $H^0_I(M)=M$.
\end{proof}

The following example  shows that in the non-local case, it is possible that  $r(I,M)=R\;$ but $\;H^0_I(M)\ne M$.

\begin{example}
Let $R$ be a ring, having a nontrivial idempotent element $e$ (i.e. $e^2=e$ and $e\ne 0,1$). Let $I=\left\langle e\right\rangle R$ and $M=I$. Notice that $I^nM=I^{n+1}M$ for all $n\ge 1.$ Therefore, by  Proposition \ref {basicreq}, we have $r(I,M)=R$. Note that $e$ cannot be killed by any power of $I$ as $I^n\cdot e=I^{n+1}=I\ne 0$. So $e \notin {H^0_I(M)}$. But $e\in M$. Hence $M\ne {H^0_I(M)}$.
\end{example}

\section{The Filtration ${\mathcal F}^I_M=\{ r(I^n,M) \}_{n\ge 0}$}
\noindent This section deals with ${\mathcal F}^I_M=\{ r(I^n,M) \}_{n\ge 0}$. We first show that it is a filtration of ideals and also an $I$-filtration. We explore its relation  with the Ratiff-Rush filtration of $M$ with respect to $I$. We also prove that if $\grade(I,M)>0$ then $r(I,M)=({\widetilde I}M:M)$.
 
For the definition of filtration of ideals, see \cite[4.4]{BH}. The following theorem shows that the collection ${\mathcal F}^I_M=\{ r(I^n,M) \}_{n\ge 0}$ of ideals is an $I$-filtration.

\begin{theorem}\label{1}
 For any $R$-module $M$, the sequence ${\mathcal F}^I_M$ is a filtration of ideals. It is also an $I$-filtration. 
\end{theorem}
\begin{proof}
 It is easy to show that 
$$r(I^n,M)\supseteq r(I^{n+1},M)\quad \mbox{for all}\; n\ge 0.$$
Next we prove the following 
\begin{itemize}
\item [$($\rm a$)$] $r(I^n,M)\cdot r(I^m,M) \subseteq r(I^{n+m},M)$ \quad for all $n,m\ge 1$.
\item [$($\rm b$)$] $I\cdot r(I^n,M) \subseteq  r(I^{n+1},M)$ \quad for all $n\ge 1$.
\end{itemize} 
(a): Let $x\in r(I^n,M)$ and $y\in r(I^m,M)$. We have
\begin{eqnarray*}
x{I^k}M &\subseteq & I^{n+k}M \quad \mbox{for all}\;k\gg 0,\\
\mbox{and}\;\;y{I^k}M  &\subseteq & I^{m+k}M \quad \mbox{for all}\;k\gg 0.
\end{eqnarray*}
Therefore, for all $k\gg 0$, we have
\begin{eqnarray*}
 xy{I^k}M &\subseteq & xI^{m+k}M \subseteq I^{n+m+k}M.\\
\mbox{Thus,}\;\;\;\;\;\;\;xy &\in & r(I^{n+m},M).
\end{eqnarray*} 
(b): This follows from (a) since $ I \subseteq r(I,M) $.
Thus ${\mathcal F}^I_M$ is an $I$-filtration.
\end{proof}

Let us recall the definition (see \cite{HJLK}, also see \cite{Nagh}) of the Ratliff-Rush submodule of $M$ with respect to $I$.

\begin{definition}  Consider the following chain of submodules of $M$ :
$$ IM\subseteq (I^2M:_M I)\subseteq (I^3M:_M I^2)\subseteq ...\subseteq (I^{n+1}M:_M I^n)\subseteq ...$$
Since $M$ is Noetherian, this chain of submodules stabilizes. We denote its stable value by $\widetilde{IM}$. We call $\widetilde{IM}$ to be the {\it Ratliff-Rush submodule of M associated with I}. The filtration $\{\widetilde{I^nM}\}_{n\ge 1}$ is called the {\it Ratliff-Rush filtration of $M$ with respect to $I$}.   
\end{definition}

\noindent {\bf Notation:} To facilitate further calculations, set
\begin{eqnarray*}
{\mathcal F}_M^I=\{r(I^n,M)\}_{n\ge 0}\;\mbox{,}\;\;\;&&{\mathcal R}({\mathcal F}^I_M) = \bigoplus_{n\ge 0} r(I^n,M)  \;\mbox{,}\;\;\;\;\;\; \\
{\mathcal R}(I)=\bigoplus_{n\ge 0} I^n\;\mbox{,}\;\;\;\;\;\;\;
G_I(R) = \bigoplus_{n\ge 0} & I^n/{I^{n+1}}& \mbox{,}\;\;\;\;\; G_I(R)_+ = \bigoplus_{n\ge 1} I^n/{I^{n+1}},\\
\;\;\;\; \;\;\;\;\;G_{I}(M)=\bigoplus_{n\ge 0}I^nM/{I^{n+1}M}\;\;\;\;\;&\mbox{and}&\;\;\;\; {\widetilde{\mathcal R}}(I,M)=\bigoplus_{n\ge 0} {\widetilde {I^nM}}.
\end{eqnarray*}

\begin{remark}
By Theorem \ref{1}, ${\mathcal F}^I_M$ is a filtration of ideals in $R$, so ${\mathcal R}({\mathcal F}^I_M)$ is a ring. Clearly $\mathcal R(I)\subseteq {\mathcal R}({\mathcal F}^I_M)$ is a subring. Since ${\mathcal F}^I_M$ is an $I$-filtration then ${\mathcal R}({\mathcal F}^I_M)$ is an $\mathcal R(I)$- module.
\end{remark}

We study the relation between the $\mathcal R(I)$-algebra ${\mathcal R}({\mathcal F}^I_M)$ and the $\mathcal R(I)$-module ${\widetilde{\mathcal R}}(M)$. We first show that

\begin{proposition}\label{6}
${\widetilde {\mathcal R}}(I,M)$ is a graded  $\;{\mathcal R}({\mathcal F}^I_M)$-module. 
\end{proposition}
\begin{proof}
 Set $I_n=r(I^n,M)$. It is enough to check that 
\begin{equation*}
I_n\cdot {\widetilde {I^mM}} \subseteq {\widetilde {I^{n+m}M}}\;\;\mbox{for all}\;\;n,m.\tag{*}
\end{equation*}
Take any $x\in I_n$ and $z\in {\widetilde {I^mM}}$. Then, by \cite[2.2(iii)]{Nagh}    and definition of $I_n$, we have
\begin{eqnarray*}
I^{k}z \subseteq  I^{k+m}M \;&\mbox{and}&\;xI^{k}M \subseteq  I^{k+n}M \quad \mbox{for all}           \;k\gg 0.
\end{eqnarray*}
Therefore, for all $k\gg 0$ 
\begin{eqnarray*}
xI^kz &\subseteq & x\cdot I^{m+k}M\subseteq   I^{n+m+k}M.\\
\mbox{So}\;\;\;xz &\in & {\widetilde {I^{n+m}M}}.
\end{eqnarray*}
\end{proof}

\begin{corollary}\label{27}
Set $J=r(I,M)$. Then 
$$I^nM \subseteq J^nM\subseteq \widetilde{I^nM}\quad for\; all\;n\ge 1.$$
\end{corollary}
\begin{proof} The first inclusion is clear. The second inclusion follows from the fact $J^n\subseteq r(I^n,M)$ and Proposition \ref{6}.
\end{proof}

\begin{corollary}\label{25}
 Set $J=r(I,M)$. Assume that $\;\grade(I,M)>0$. Then 
$$J^nM=I^nM\quad\mbox{for all}\;n\gg 0.$$
Furthermore if $\grade\left(G_I(R)_{+},G_I(M)\right) >0$ then
$$J^nM=I^nM\quad\mbox{for all}\; n\ge 1,$$
\end{corollary}
\begin{proof}
 First part follows from \cite[2.2]{Nagh}. For second part, note that
$\grade(G_I(R)_{+},$\;\\$G_I(M))>0\;$ implies \;(see \cite[Fact 9]{HJLK}) that
$\widetilde{I^nM}=I^nM\quad\mbox{for all}\;n\ge 1.$ 
\end{proof}

\begin{proposition}\label{relationf}
 If $\;\grade(I,M)>0$ then $$r(I,M)= (\widetilde{IM} :M).$$
\end{proposition}
\begin{proof} From Proposition \ref{6}, it follows that $r(I,M)\subseteq (\widetilde{IM} :M)$. For reverse inclusion, let $x\in (\widetilde{IM} :M)$. Then
\begin{eqnarray*}
xM &\subseteq & \widetilde{IM} ,\\
\mbox{so}\quad xI^kM &\subseteq & I^k \widetilde{IM} = I^{k+1}M\quad \mbox{for all} \;\;k\gg 0.
\end{eqnarray*}
Thus $\; x \in r(I,M)$. Hence $r(I,M)= (\widetilde{IM} :M).$ 
\end{proof}

\section{Involution Properties}
\noindent In this section we prove that the function $\;I\longmapsto r(I,M)\;$ is an involution on the set of ideals of $R$ (see Theorem \ref{closedness2}). We first prove the result in the case when $\grade(I,M)>0$. We also show that if $\;\grade(I,M)>0$ then $\;r(I,M)$ is a Ratliff-Rush closed ideal.
   
\begin{proposition}\label{closedness1} 
Let $\grade(I,M)>0$. Set $J=r(I,M)$. Then 
\begin{itemize}
\item [$($\rm a$)$] $r\left(J,M \right) =  J.$
\item [$($\rm b$)$] $\widetilde{J}= J.$
\end{itemize}
\end{proposition}
\begin{proof}
 {\bf (a)} By Corollary \ref{25}, there exists an integer $k_0$ such that
\begin{equation*}
J^kM=I^kM\quad\mbox{for}\;k\ge k_0.\tag{i}
\end{equation*}
Also there exists ${k_0}'$ such that
\begin{eqnarray*}
r(J,M) &=& (J^{k+1}M: J^kM)\quad\mbox{for all}\;k\ge {k_0}'.
\end{eqnarray*}
From (i), it follows that 
\begin{eqnarray*}
r(J,M)&=& (I^{k+1}M: I^kM) \quad \mbox{for all}\;k\ge \max\{k_0,{k_0}'\}\\
&\subseteq & r(I,M)\;=\;J.
\end{eqnarray*}
Since $J\subseteq r(J,M)$ always we get $r(J,M)=J$.\\
{\bf (b)} Let $x\in \widetilde J$. Then
\begin{eqnarray*}
xJ^k &\subseteq & J^{k+1}\;\;\mbox{for all}\;k\gg 0,\\
xJ^kM &\subseteq &  J^{k+1}M \;\;\mbox{for all}\;k\gg 0.
\end{eqnarray*}
By Corollary \ref{25}, we get
\begin{eqnarray*}
 xI^kM &\subseteq &  I^{k+1}M \;\;\mbox{for all}\;k\gg 0.
\end{eqnarray*}
Therefore $x\in J$. Hence $\widetilde J=J$.
\end{proof}

To tackle the case when $\grade(I,M)=0$ we first prove

\begin{lemma}\label{30}
For any ideal $I$, the following hold
\begin{itemize}
\item [$($\rm a$)$] $r(I,M) = r\left(I,M/{H^0_I(M)}\right).$
\item [$($\rm b$)$] $r(I^n,M) = r\left(I^n,M/{H^0_I(M)}\right)\;\;\mbox{for all}\;\;n\ge 1.$
\end{itemize}  
\end{lemma}
\begin{proof}
{\bf (a)} If $M=H^0_I(M)$ then our assertion follows from Proposition \ref{basicreq}.
Suppose $M\ne H^0_I(M)$. Since the natural map  $M\longrightarrow M/{H^0_I(M)}$ is surjective,  by Proposition \ref{3}(a), we have 
$$r(I,M)\subseteq r\left(I,M/{H^0_I(M)}\right).$$
Let $x\in r\left(I,M/{H^0_I(M)}\right)$. Therefore
\begin{eqnarray*}
xI^kM+H^0_I(M) &\subseteq & I^{k+1}M+H^0_I(M)\quad\mbox{for some}\; k\ge 0.
\end{eqnarray*}
Since $M$ is a finitely generated $R$-module, there exists an integer $r\in \mathbb{N}$ such that
$I^rH^0_I(M)=0$. Therefore
\begin{eqnarray*}
xI^{r+k}M &\subseteq & I^{r+k+1}M.
\end{eqnarray*}
So $x \in  (I^{r+k+1}M: I^{r+k}M)\subseteq r(I,M)$.\\
\noindent{\bf (b)} This follows from (a), since
$\;H^0_{I^n}(M)=H^0_I(M)\quad\mbox{for all}\; n\ge 1.$   
\end{proof}

We now prove the involution property in general.

\begin{theorem}\label{closedness2}
 For any ideal $I$ of ring $R$, we have
$$r\left(r(I,M), M\right) = r(I,M).$$

\end{theorem}
\begin{proof} Set $J=r(I,M) \;\;\mbox{and}\;\; N=M/{H^0_I(M)}.$
If $M=H^0_I(M)$, then $J=R$ by Proposition \ref{basicreq}. Also clearly $r(R,M)=R$.
 Now suppose $M\ne H^0_I(M)$. By using Lemma \ref{30} for the ideals $I$ and $J$, we have
\begin{equation*}
J = r(I,M) \;= \;r(I,N).\tag{i}
\end{equation*}
\begin{equation*}
r(J, M) = r\left( J,M/{H^0_J(M)}\right)\tag{ii}.
\end{equation*}
Since $I\subseteq J$ we get $H^0_J(M)\subseteq H^0_I(M).$
Therefore the epimorphism 
$$M/{H^0_J(M)}\longrightarrow M/{H^0_I(M)}=N$$ 
induces, by Proposition \ref{3}(a),
\begin{equation*}
r(J,M/{H^0_J(M)}) \subseteq  r(J,N).\tag{iii}
\end{equation*}
Notice that $\grade(I, N)>0$. Therefore, from Proposition \ref{closedness1}, we get 
\begin{equation*}
r(J,N) = r\left( r(I,N),N\right)  = r(I,N).\tag{iv}
\end{equation*}
Using (ii), (iii) and (iv), we get
\begin{eqnarray*}
r(J,M) &=& r\left(J,M/{H^0_J(M)}\right)\;\;\subseteq\;\; r(J,N)=J.
\end{eqnarray*}
But $r(J,M) \supseteq J$ always. Therefore $r(J,M)=J$. 
\end{proof}

\section{Relation with Integral Closure}
\noindent In this section we show that $r(I,M) \subseteq \overline I$, the integral closure of $I$ when 
$\ann M = 0$ and $\grade(I,R) > 0$. In particular when
$M$ has a positive rank and $I$ is a regular ideal. We prove that if $I$ is a regular ideal then there exists an $R$-module $M$ of rank 1 such that $r(I,M)=\overline I$. Finally we show that if $I$ is a regular ideal then the set 
$${\mathcal C}(I):=\{\;J:\; J\mbox{ is regular ideal and}\;r(I,J)={\overline I}\;\}$$
is non-empty and  has a unique maximal element.

\begin{remark} For any $x\in r(I,M)$, there exists $k\ge 0$ such that $xI^kM\subseteq I^{k+1}M$. By determinant trick, there exists  $f(t)\in R\left[ t \right]$ such that 
$$f(t)= t^{n}+a_1t^{n-1}+...+a_{n-1}t+a_{n},\;\;\mbox{where}\;a_i\in I^i,$$
and $f(x) \in  \ann (I^kM).$
\end{remark}  
The following proposition gives a relation between $r(I,M)$ and $\overline I$.

\begin{proposition}\label{integral1}
Let $\;\ann M=0$. If either of following two conditions holds
\begin{itemize}
\item [$($\rm a$)$] $\grade(I,R)>0$. 
\item [$($\rm b$)$] $\grade(I,M)>0$.
\end{itemize}
Then $r(I, M)\subseteq \overline{I}.$ In particular if $I$ is a regular ideal and 
$M$ has a positive rank then $r(I, M)\subseteq \overline{I}.$
\end{proposition}
\begin{proof}
{\bf (a)}
By the remark above we get that $f(x)I^kM = 0$ for some $k \geq 0$. But $\ann M = 0$. So 
  $f(x)\cdot I^k = 0$.  Since $\grade(I, R)>0$, we get $f(x)=0$. Hence $x\in {\overline I}$.

\noindent {\bf (b)} Note that $\grade(I,M)>0$ yields $\ann(I^kM)=\ann M$. Now from hypothesis $\ann M=0$, it follows that $f(x)=0$. Hence $x\in {\overline I}$.

Finally note that if $M$ has a positive rank, say $r$, then $\ann M = 0$ (since it contains
$R^r$ as a submodule.
\end{proof}

\begin{remark}
 Notice that $\;\ann M=0\;$ together with $\;\grade(I,M)>0\;$ implies \\
$\;\grade(I,R) > 0$. Thus {\bf (b)} follows from {\bf (a)}. 
\end{remark}

The next proposition ensures the existence of an $R$-module $M$ for a regular ideal $I$ of $R$ such that $r(I,M)={\overline I}.$

\begin{theorem}\label{10}
Let $I$ be a regular ideal of $R$. Let $J$ be an ideal such that $I\subseteq J\subseteq \overline I$. Then there exists $R$-module $M$ of rank 1  such that $$I \subseteq J \subseteq r(I,M)\subseteq {\overline I}.$$
In particular, there exists an $R$-module $M$ of rank 1 such that $r(I,M)={\overline I}.$
\end{theorem}
\begin{proof}
Let $z\in J\setminus I$ be any element. We have $z^n+\sum_{i=1}^n a_iz^{n-i} = 0,\;\;\mbox{with}\; a_i \in I^i$. Thus
\begin{equation*}
\mbox{}\;\; z^n = - \sum_{i=1}^n a_iz^{n-i}.\tag{*}
\end{equation*}
Set $N=\left\langle z,I\right\rangle ^{n-1}R$. We claim that $zN\subseteq IN$. Let $x\in N$. Then 
\begin{eqnarray*}
x &=& rz^{n-1}+ \sum_{i=1}^{n-1} b_iz^{n-1-i},\;\;\mbox{where}\;r\in R\;\mbox{and}\;b_i\in I^i.\\
\mbox{Therefore}\;\;\;\;zx &=& rz^n+ \sum_{i=1}^{n-1} b_iz^{n-i}.
\end{eqnarray*}
By using (*), we get $zx = \sum_{i=1}^{n-1}(b_i-ra_i)z^{n-i} - ra_n \in IN$.
Thus $z N\subseteq IN$ and hence $z\in r(I,N).$\\
\indent Let us assume that $J =\left\langle z_1,z_2,...,z_s\right\rangle R$\quad and set\quad $N_i=\left\langle z_i, I\right\rangle ^{n_i-1}R$, where $n_i$ is the degree of an integral equation satisfied by $z_i$. Set $M=N_1\otimes_R N_2\otimes_R ...\otimes_R N_s$. Notice that $\rank M=1$. By Proposition \ref{3}, we get $J \subseteq r(I,M)$. But, by Proposition \ref{integral1}(a), we always have $\overline I \supseteq r(I,M)$. 
\end{proof}

\begin{proposition}\label{rank2}
For any regular ideal $I$, there exists a  regular ideal $J$ such that 
$$r(I,J)=\overline I.$$
\end{proposition}
\begin{proof}
 By Theorem \ref{10}, there exists an $R$-module $M$ of rank 1 such that  $r(I,M)=\overline I$.  Let $T(M)$ denote the torsion submodule of $M$ and 
set $N=M/{T(M)}$. Note that the surjective map $M\longrightarrow N$ induces $r(I,N) \supseteq r(I,M)=\overline I$. But as $\rank N=\rank M=1$, we have $r(I,N) \subseteq \overline I$  $\left( \mbox{by} \;\ref{integral1}(a)\right)$. 
So $r(I,N)=\overline I$. Clearly $N$ is torsion-free so that $N\simeq J$\; for some ideal $J$ of $R$. Since $J$ is of rank $1$, $J$ is a regular ideal. Notice that $r(I,J)=r(I,N)=\overline I. $ 
\end{proof}

Let ${\mathcal C}(I):=\{\;J:\; J\mbox{ is a regular ideal and}\;r(I,J)={\overline I}\;\}.$ By Proposition \ref{rank2}, ${\mathcal C}(I)\ne \emptyset.$  Since $R$ is Noetherian, ${\mathcal C}(I)$ has a maximal element. We show that  

\begin{theorem}\label{existence1}
 ${\mathcal C}(I)$ has a unique maximal  element.
\end{theorem}
\begin{proof}
 Suppose $Q\in {\mathcal C}(I)$ is a maximal element and $J\in {\mathcal C}(I)$. By Proposition \ref{3}, the following epimorphism \;$J\oplus Q\longrightarrow J+ Q\longrightarrow 0$\; induces
\begin{eqnarray*}
r(I,J\oplus Q)&\subseteq & r(I,J+ Q).\\
\mbox{But} \;\;\;r(I,J\oplus Q) &=& r(I,J)\bigcap r(I,Q)={\overline I},\quad \mbox{by Proposition \ref{3}(b).}\\
\mbox{So}\quad {\overline I} &\subseteq & r(I,J+ Q).
\end{eqnarray*}
But $r(I,J+ Q) \subseteq {\overline I}$, since $J+Q$ has rank 1 as an $R$-module so  $\;r(I,J+Q)={\overline I}$. As $Q$ is maximal this gives $Q=J+Q$. So $J\subseteq Q$.  Hence ${\mathcal C}(I)$ has a unique maximal  element.
\end{proof}

\section{Stable filtrations}
 In this section  we  discuss  the conditions under which our filtration $\mathcal{F}^I_M =\{ r(I^n,M) \}_{n\ge 0}$ is a stable $I$-filtration. This is equivalent to saying that the Rees algebra ${\mathcal R}({\mathcal F}^I_M)$ is a finitely generated ${\mathcal R}(I)$-module. Our main result (Theorem \ref{stablity2}) is that if $\grade(I,R)>0$ and $\ann M=0$ then ${\mathcal F}^I_M$ is a stable $I$-filtration. In local case we  prove that $\ann M=0$ is a necessary condition for ${\mathcal F}^I_M$ to be a stable $I$-filtration.

\s Recall  a filtration of ideals
$R=I_0 \supseteq  I_1 \supseteq  ...\supseteq I_n \supseteq I_{n+1} \supseteq... $ is said to be a {\it stable $I$-filtration} if $\;II_n\subseteq I_{n+1}$ for all $n\ge 0$ and $\;II_n= I_{n+1}$ for $n\gg 0$.

The lemma below is crucial to prove our main result.
\begin{lemma}\label{finiteness1}
Let $S$ be a ring and $R\subseteq S$, a subring of $S$, such that $R$ is Noetherian. Assume that there is a faithful $S$-module $E$ $($i.e, $\ann_S(E)=0$ $)$ such that $E$ is a finitely generated  $R$-module. Then $S$ is finitely generated as a $R$-module (and so Noetherian).
\end{lemma}
\begin{proof}
Note that any $S$-linear map  $f:M\longrightarrow N$  between $S$-modules $M$ and $N$, is  also $R$-linear. Consider the inclusion map
$$i:{\Hom}_S(E,E)\longrightarrow {\Hom}_R(E,E)\;\;\mbox{such that}\;\; f\longmapsto f.$$
Notice that $i$ is $R$-linear. For  $s\in S$, let $\mu_s :E\longrightarrow E$ be the multiplication map i.e., $\mu_s(t)=st\;\;\mbox{for all}\; t\in E$. Define 
\begin{eqnarray*}
\phi: S &\longrightarrow & {\Hom}_S(E,E)\\
\;\;\; s &\longmapsto & \mu_s .
\end{eqnarray*}
Clearly $\phi$ is $S$-linear and so $R$-linear. Notice that $\ker{\phi} =0$, since $E$ is a faithful $S$-module.
Consider the following composition 
\begin{equation*}
S \xrightarrow{\phi} {\Hom}_S(E,E) \xrightarrow{i} {\Hom}_R(E,E).
\end{equation*}
Clearly $i\circ \phi$ is an injective $R$-linear map. Therefore as $R$-modules
$$ S\cong \ \text{to a} \  R\text{-submodule of} \  \Hom_R(E,E).$$ 
As $R$ is Noetherian and $E$ is a finitely generated $R$-module, we get $S$ is a finitely generated $R$-module. 
\end{proof}

The next theorem shows that the filtration ${\mathcal F}^I_M$ is a stable $I$-filtration under fairly mild assumptions.

\begin{theorem}\label{stablity1}
Let $\grade(I,M)>0$ and $\;\ann M=0$. Then ${\mathcal F}^I_M$ is a stable $I$-filtration.
\end{theorem}
\begin{proof}
 For convenience, set $\mathcal S={\mathcal R}({\mathcal F}^I_M)$, $\;{\mathcal{R}}={\mathcal R}(I)$ and $E=\oplus_{n\ge 0} {\widetilde {I^nM}}t^n$. By Proposition \ref{6}, $E$ is an $\mathcal S$-module. Since $\grade(I,M)>0$,  ${\widetilde {I^nM}}=I^nM$ for all $n\gg 0$ (see \cite[3.3]{Nagh}). So $E$ is a finitely generated ${\mathcal R}$-module. 

We  prove that $\ann_{\mathcal S}(E)=0$. Notice that $\ann_{\mathcal S}(E)$ is a homogeneous ideal of $\mathcal S$.  Let $xt^n\in \ann_{\mathcal S}(E)$ be a homogeneous element. As $xt^n\cdot E=0$ \;we get  $x\cdot M=0$. Thus $x\in \ann M=0$. Therefore $\ann_{\mathcal S}(E)=0$. Using Lemma \ref{finiteness1}, we conclude that $\mathcal S$ is a finitely generated ${\mathcal R}$-module. So  ${\mathcal F}^I_M$ is a stable $I$-filtration.
\end{proof}

The following example shows that the hypothesis in Theorem \ref{stablity1} is not necessary for ${\mathcal R}({\mathcal F}^I_M)$ to be Noetherian.

\begin{example} Let $I$ be a nilpotent ideal of $R$ (i.e. $I^r=0$ for some $r\ge 1$). Then $r(I^n,M)=R$ for all $n\ge 0$ and for any $R$-module $M$. Thus ${\mathcal R}({\mathcal F}^I_M) \cong R[t]$, which is a Noetherian ring.
\end{example}

\begin{remark}\label{annconinrat}
Let $I$ be any ideal of $R$. Let $x\in \ann M$. Then $x\cdot M=0$, so $xI^kM=0\subseteq I^{n+k}M$ for all $n,k\ge 1$. Therefore $x\in r(I^n,M)$ for all $n\ge 1$. Hence 
$$\ann M\subseteq \bigcap_{n\ge 1}r(I^n,M).$$
\end{remark}

\begin{proposition}\label{annconinter}
If ${\mathcal F}^I_M$ is a stable $I$-filtration then
$$\ann M\subseteq \bigcap_{n\ge 1}I^n.$$
\end{proposition}
\begin{proof}
 For convenience, set $I_n=r(I^n,M)$.
Since ${\mathcal F}^I_M$ is a stable $I$-filtration, there exists an integer $n_0$ such that
\begin{equation*}
I_{n_0+k}=I^kI_{n_0} \quad \mbox{for all}\; k\ge 1.\tag{i}
\end{equation*}
By above Remark \ref{annconinrat}, one has $\ann M\subseteq I_{n_0+k}$ for all $k\ge 1$. Therefore, by (i), we have
$$\ann M\subseteq I_{n_0+k} \subseteq I^k\quad\mbox{for all}\;k\ge 1.$$
Hence result follows.
\end{proof}

\begin{remark} 
The above Proposition \ref{annconinter} proves that if the $I$-adic filtration  is {\it separated}, i.e. $\bigcap_{n\ge 1}I^n=0$ then 
\[
{\mathcal F}^I_M\;\;is \;\;a \;stable \;I\mbox{-filtration}\;\; \Longrightarrow \;\;\ann M=0.
\]
\end{remark}

An easy consequence of Proposition \ref{annconinter} is following

\begin{corollary}\label{Noeth}
Let $(R, \m)$ be a local ring and $M\ne H^0_I(M)$. Then
$${\mathcal F}^I_M\;\;is \;\;a\;stable\;I\mbox{-filtration}\;\; \Longrightarrow \;\;\ann M=0.$$
\end{corollary}
\begin{proof}
 For a local ring $(R,\m)$, the $I$-adic filtration is separated (by Krull's intersection theorem). 
Hence $\ann M=0.$ 
\end{proof}

In the next proposition we prove a partial converse of above Corollary \ref{Noeth}. 
\begin{theorem}\label{stablity2} 
Let $I$ be a regular ideal. If $\;\ann M=0$ then ${\mathcal F}^I_M$ is a stable $I$-filtration.
\end{theorem}
\begin{proof}
 Notice that $M\ne {H^0_I(M)}$. Set $N=M/{H^0_I(M)}$. Note that $\grade(I,N)>0$. We have  
\begin{equation*}
r(I^n,M) = r\left(I^n,N\right)\quad\mbox{for all}\;n\ge 0.\quad \mbox{(by Lemma \ref{30})} \tag{*}
\end{equation*} 
Let $x\in \ann_R(N)$. We have $xM \subseteq H^0_I(M)$. Thus there exists $k\ge 1$ such that $I^k(xM)=0$, so $xI^k\subseteq \ann M=0$. But since $I$ is regular, $x=0$. Hence $\ann_R(N)=0$.  Therefore, by Theorem \ref{stablity1}, the filtration $\{r(I^n,N)\}_{n\ge 0}$ is a stable $I$-filtration and so is ${\mathcal F}^I_M$.
\end{proof}

\section{The case when $M_{\p}$ is free for all $\p \in \spec(R) \setminus \mspec(R)$.}
\noindent In this section we study our filtration ${\mathcal F}^I_M=\{r(I^n,M)\}_{n\geqslant 1}$ when $M$ is free for all $\p\in \spec R\setminus \mspec R.$ We show that for a regular ideal $I$, ${\mathcal F}^I_M$ is a stable $I$-filtration when $\Ass R\cap \mspec R=\emptyset$. We also prove that if $ A^*(I)\cap \mspec R=\emptyset\;$ then $\;r(I^n,M)=\widetilde{I^n}$\;for all\; $n\geqslant 1.$ Here $A^*(I)$ is the stable value of the sequence $\Ass(A/{I^n})$.

\s \label{asympassump}
Throughout this section we assume that
\begin{enumerate}
\item [$($\rm 1$)$] $\Ass R\cap \mspec R=\emptyset$ and 
\item [$($\rm 2$)$] $M$ is an $R$-module such that $M_{\p}$ is free for all $\p \in \spec R\setminus \mspec R$, where $\mspec R=\{ \m\;:\;\m \;\mbox{is a maximal ideal of}\; R\}.$
\end{enumerate}
\indent We give some examples where these assumptions hold.

\begin{examples}$~~$
\begin{itemize}
\item [$($\rm 1$)$] $\Ass R\cap \mspec R=\emptyset$ holds if and only if $\depth R_{\m}>0$ for all $\m \in \mspec R$. Thus if $(R,\m)$ is a local domain which is not a field then assumption \ref{asympassump}(1) holds.
\item [$($\rm 2$)$] If an $R$-module $M$ satisfies the exact sequence of the form
$$0\longrightarrow M \longrightarrow F \longrightarrow F/M\longrightarrow 0, \quad \mbox{with}\;\;\ell(F/M)<\infty,$$
where $F$ is a free $R$-module then $M_{\p}$ is free for all $\p \in \spec R\setminus \mspec R$.
\item [$($\rm 3$)$] Let $(R,\m)$ be a local Cohen-Macaulay ring and an isolated singularity i.e.,  $R_{\p}$ is regular local for all prime $\p\ne \m.$
Then if $M$ is a maximal Cohen-Macaulay $R$-module then $M_{\p}$ is free for all $\p\ne \m.$
\end{itemize}
\end{examples}

\begin{lemma}\label{asymlem1}
$($with  hypotheses as in \ref{asympassump}$)$ $\;\ann M=0$.
\end{lemma}
\begin{proof} Notice that $(\ann M)_{\p}=\ann_{R_{\p}}(M_{\p})=0$, for all $\p \in \spec R\setminus \mspec R.$ By our first hypothesis, we can have $\Ass (\ann M)\cap \mspec R=\emptyset$ and so $\Ass (\ann M)=\emptyset$. Therefore $\ann M=0$.
\end{proof}
The following proposition readily follows from Theorem \ref{stablity2} and  Lemma \ref{asymlem1}.

\begin{proposition}\label{stabilitysecond}
$($with  hypotheses as in \ref{asympassump}$)$ If $I$ is a regular ideal then ${\mathcal F}^I_M$ is a stable $I$-filtration and so $\mathcal{R}({\mathcal F}^I_M)$ is  finitely generated as an ${\mathcal R}(I)$-module.
\end{proposition}

\begin{remarks}\label{asymrem1}$~~$
\begin{itemize}
\item [$($\rm 1$)$] By the result of Brodmann \cite{Brod}, the sequence $\Ass(R/{I^n})$ stabilizes for large $n$. Let $A^*(I)$ denote the stable value of this sequence.
\item [$($\rm 2$)$] Ratliff in his paper \cite[2.7]{Rat}, has proved that the sequence $\Ass (R/{\overline {I^n}})$ eventually stabilizes at a set denoted by ${\overline {A^*}}(I)$.
\item [$($\rm 3$)$] In \cite[2.8]{Rat}, it is also proved that $\overline {A^*}(I) \subseteq A^*(I)$.
\item [$($\rm 4$)$] By \cite[1.6]{Mc}, we have $\p \in A^*(I)$ if and only if $\p_S\in A^*(I_S)$, for any multiplication closed set $S$ disjoint from $\p$. 
\end{itemize}
\end{remarks}
The following is well-known. We include a proof for lack of a suitable reference.
\begin{lemma}\label{asymplem2}
For a regular ideal $I$ we have
\begin{itemize}
\item [$($\rm a$)$] $\Ass(R/{\widetilde{I^n}})\subseteq \Ass(R/{\widetilde{I^{n+1}}})$\; for all $n\geqslant 1.$
\item [$($\rm b$)$] Further, $\Ass(R/{\widetilde{I^n}})\subseteq A^*(I)$ \; for all $n\geqslant 1.$
\end{itemize}
\end{lemma}
\begin{proof} {\bf (a)} Fix $n\geqslant 1$. Let $\p \supseteq I$ be such that $\p \in \Ass(R/{\widetilde{I^n}})$. 
We localize $R$ at $\p$. Set $\m=\p R_{\p}$. Since associated primes behave well with respect to localization so we may assume that $(R,\m)$ is local and $\m\in \Ass(R/{\widetilde{I^n}})$. We may further assume that $R/{\m}$ is infinite. Otherwise we make a base change $R\longrightarrow R[X]_{\m R[X]}=T$. Let $\n= \m T$, the extension of the maximal ideal of $R$ in $T$. Notice that if $E$ is an $R$-module then 
$$\n\in \Ass_T(E\otimes_R T) \;\;\mbox{if and only if}\;\;  \m\in \Ass E.$$
Therefore we assume that $(R,\m)$ is local with $R/{\m}$ infinite and $\m\in \Ass(R/{\widetilde{I^n}})$.
 Let $x\in I$ be a superficial element with respect to $I$. Consider the map
$$\mu^x_n: R/{\widetilde{I^n}} \longrightarrow R/{\widetilde{I^{n+1}}},\quad \mbox{such that}\;\;a+\widetilde{I^n} \longmapsto ax+\widetilde{I^{n+1}}.$$
Clearly $\mu^x_n$ is $R$-linear. Also it is  injective.  Thus $\;\m\in \Ass(R/{\widetilde{I^n}})\subseteq \Ass(R/{\widetilde{I^{n+1}}})$.\\
\noindent {\bf (b)} By repeatedly using (a) we get
$$\p\in \Ass(R/{\widetilde{I^{n+k}}})\;\;\mbox{for all}\;\;k\ge 1.$$
Note that for $k\gg 0$,  $\widetilde{I^{n+k}}=I^{n+k}$. Also by Remark \ref{asymrem1}(1), $\Ass(R/{I^{n+k}})=A^*(I)$. Therefore the result follows.
\end{proof}

\begin{theorem}\label{asymptheo1}
$($with  hypotheses as in \ref{asympassump}$)$ Let $I$ be a regular ideal of $R$. Then the function $n\mapsto \ell ( r(I^n,M)/{\widetilde{I^n}})$ is a polynomial function.  Furthermore if $  A^*(I)\cap \mspec R=\emptyset$ then $r(I^n,M)=\widetilde{I^n}$ for all $ n\geqslant 1.$
\end{theorem}
\begin{proof}
Notice that $\widetilde {I^n}\subseteq r(I^n,M) \ \mbox{for all} \  n\ge 1.$
Let $\p\in \spec R\setminus \mspec R$. By hypotheses, Proposition \ref{3}(e) and Remark \ref{prelifreemodule}, we have 
\begin{equation*}
r(I^n,M)_{\p}=r_{R_{\p}}(I^n_{\p},M_{\p})=\widetilde{I^n_{\p}}\;\;\mbox{for all}\;n\ge 1. \tag{*}
\end{equation*}
Therefore for all $n\geq 1$, 
$$\ell\left( \frac{r(I^n,M)}{\widetilde {I^n}} \right) \quad \mbox{is finite}.$$ 
By Proposition \ref{stabilitysecond}, $\mathcal{R}({\mathcal F}^I_M)$ is a finitely generated ${\mathcal R}(I)$-module. Also  $\widetilde{\mathcal R}(I)$ is a finitely generated $\mathcal R(I)$-module. So
$$E = \frac{\mathcal{R}({\mathcal F}^I_M)}{\widetilde{\mathcal R}(I)} = \bigoplus_{n\geq 1}\frac{r(I^n,M)}{\widetilde {I^n}} $$
is  a finitely generated $\mathcal{R}(I)$-module. From a well-known fact it follows  that the function
$n\longmapsto \ell\left( r(I^n,M)/{\widetilde {I^n}}\right)$ is a polynomial function.\\
\indent The exact sequence 
$0\rightarrow r(I^n,M)/{\widetilde {I^n}}\rightarrow R/{\widetilde {I^n}}$
yields $$\Ass\left( \frac{r(I^n,M)}{\widetilde {I^n}} \right) \subseteq \Ass \left( R/{\widetilde {I^n}}\right) \subseteq A^*(I).$$
If $ A^*(I)\cap \mspec R=\emptyset$ then by (*) we get $r(I^n,M)=\widetilde{I^n}\;$\;for all\; $n\geqslant 1.$
\end{proof}

\begin{remark}\label{rrneqM}
 If  $A^*(I)\cap \mspec R \neq \emptyset $ then $r(I,M)$ need not be equal to $\widetilde{I}$. For instance
 let $R = k[t^4,t^{11},t^{17}, t^{18}]_{\m}$ where $\m = ( t^4,t^{11},t^{17}, t^{18})$. Let $I = (t^4,t^{11})$. Set $M = (t^4,t^{11},t^{17})$ considered as a submodule
of $R$. Notice that $\ell (R/M)$ is finite. By Example \ref{examplerepeat1}  we get that $r(I^n,M) \neq \widetilde{I^n}$ for all $n \geq 1$. In this case
 $I$ is $\m$-primary. So  $A^*(I) = \{ \m \}$.
 \end{remark}

In view of Theorem \ref{asymptheo1} we give some situations of prime $\p$ such that $\p \in A^*(I)$.

\begin{remarks}\label{asymplastremark}$~~$
\begin{itemize}
\item [$($\rm 1$)$] If $ht(\p)=l(I_{\p}),$  the \emph{analytic spread} of $I_{\p}$ then $\p \in A^*(I).$ (see \cite[4.1]{Mc})
\item [$($\rm 2$)$] If $\p \in A^*(I) \setminus B^*(I)$, where $B^*(I)$ is the stable value of the sequence $\Ass (I^n/{I^{n+1}})$, then $\p \in \Ass R$. (see \cite[2.2]{Mc})
\item [$($\rm 3$)$] With  hypotheses as in \ref{asympassump}, if a maximal ideal $\m\in A^*(I)$ then $\m \in B^*(I)$ that is, $\m\in \Ass_R(I^n/{I^{n+1}})$ for $n\gg 0$. So we get $\m/I \in \Ass_{R/I}(I^n/{I^{n+1}})$ for $n\gg 0.$ Set $G = G_I(R)$. By \cite[2.1]{Sche}, we thus have $\m/I = Q\cap R/I$ such that $Q \in \Ass\left( G\right)\setminus V( G_+).$
\end{itemize}
\end{remarks}

\begin{corollary}
$($with  hypotheses as in \ref{asympassump}$)$  In addition let $(R,\m)$ be a local Cohen-Macaulay ring. Let $x_1,\ldots,x_r$ be a regular sequence with $r<\dim R$. Set $I=\left\langle x_1,\ldots,x_r \right\rangle R$. Then $r(I^n,M)=I^n$ for all $n\geq 1$.
\end{corollary}
\begin{proof}
Clearly $R/{I^n}$ is a Cohen-Macaulay ring of dimension greater than or equal to $1$ for all $n\geqslant 1$. So $\m\notin \Ass(R/{I^n})$ for all $n\geqslant 1$. This gives $\m\notin A^*(I)$. Therefore by Theorem \ref{asymptheo1}, 
$r(I^n,M)= \widetilde{I^n}$ for all $n\geq 1$. However as $\depth G_I(R) > 0$ we get $\widetilde{I^n} = I^n$ for all
$n \geq 1$. The result follows.
\end{proof}

\section{Some more Analysis on $\;r(I,M)$}
\noindent In this section we analyze the case when $\ann M$ need not be zero. We also consider the case when $\grade(I,M)=0$. When $M\ne H^0_I(M)$ both these cases can be dealt with by going modulo the ideal $(H^0_I(M): M)$. We prove that $r(I^{n+1},M)=I\cdot r(I^n,M)+(H^0_I(M):M)$\;\;for all \;$n\gg 0.$
Our techniques also yield $\widetilde{I^nM}=I^nM+H^0_I(M)\;\;\mbox{for all}\;n\gg 0.$

\s Before we proceed further let us fix some notations which we will use throughout the section. Set $N=M/{H^0_I(M)}$, $q_{_I}(M)=(H^0_I(M): M)$ and $S=R/{q_{_I}(M)}$. Let $J$ be the image of $I$ in $S$.
 
\begin{proposition}\label{recursion1}
Let $M$ be an $R$-module such that $M\ne H^0_I(M)$. Then 
\begin{itemize}
\item [$($\rm a$)$] $r(I^n,M)=r_R(I^n,N)\;$ for all $\;n\ge 1.$
\item [$($\rm b$)$] $q_{_I}(M) \subseteq r(I^n,N)\;$ for all $\;n\ge 1.$
\item [$($\rm c$)$] $\ann_S(N)=0$.
\item [$($\rm d$)$] $\grade_R(I,N)=\grade_S(J,N)>0$.
\item [$($\rm e$)$] $r(I^n,M)/{q_{_I}(M)}=r_S(J^n,N)\;$ for all $\;n\ge 1.$
\end{itemize}
\end{proposition}
\begin{proof} Set $q=q_{_I}(M)$. (a) follows from Lemma \ref{30}(a). Parts  (b) and (c) are easy to prove. For (d), note that $H_J^0(N)=0$. To prove (e), it is sufficient to show that
$$\frac{r(I^n,N)}{q}=r_S(J^n,N)\quad\mbox{for all}\; n\ge 1.$$
Let $x\in r(I^n,N)$. Thus $xI^kN\subseteq I^{n+k}N,$ so by going modulo $q$, we get
\begin{eqnarray*}
{\bar x} J^kN &\subseteq & J^{n+k}N.\;\mbox{So}\quad {\bar x} \in  r_S(J^n,N).
\end{eqnarray*}
Conversely if ${\bar x}\in r_S(J^n,N)$ then we have ${\bar x} J^kN \subseteq  J^{n+k}N.$ Thus
\begin{eqnarray*}  
(x+q)\left( \frac{I^k+q}{q}\right)N  &\subseteq & \left( \frac{I^{n+k}+q}{q}\right)N,\\
\mbox{so}\;\;\;\left( \frac{xI^k+q}{q}\right)N  &\subseteq & \left( \frac{I^{n+k}+q}{q}\right)N.
\end{eqnarray*}
This implies $\; ({xI^k+q})N \subseteq  ({I^{n+k}+q})N.$ So $\;xI^kN \subseteq I^{n+k}N,$
since $q=\ann_R N$. Therefore $x\in r(I^n,N)$ and hence ${\bar x}\in r(I^n,N)/{q}$.
\end{proof}

\begin{theorem}\label{recursion}
Let $M$ be an $R$-module such that $M\ne H^0_I(M)$. Then 
$$r(I^{n+1},M)=I\cdot r(I^n,M)+(H^0_I(M):M)\quad \mbox{for all }\;n\gg 0.$$
\end{theorem}
\begin{proof}
By Proposition \ref{recursion1}(d), $\grade_S(J,N)>0$. Together with result (c), this gives that the filtration $\{r_S(J^n,N)\}_{n\ge 0}$ is a stable $J$-filtration (by Theorem \ref{stablity1}). Therefore we have
$r_S(J^{n+1},N) = J\cdot r_S(J^n,N)\;\;\mbox{for all}\;n\gg 0$. So
\begin{eqnarray*}
 \frac{r(I^{n+1},M)}{q} &=& \left(\frac{I+q}{q}\right) \cdot \left( \frac{r(I^n,M)}{q}\right).
\end{eqnarray*}
Thus $\;r(I^{n+1},M) = I\cdot r(I^n,M)+q \;\;\mbox{for all}\;n\gg 0$.
\end{proof}

\s {\bf Consequences of Theorem \ref{recursion}}
\begin{itemize}
\item [$($\rm 1$)$] When $M=R$ we have, for any ideal $I$,
$$\widetilde{I^{n+1}}=I\cdot \widetilde{I^n}+H^0_I(R)\;\;\mbox{for all}\;n\gg 0.$$
\item [$($\rm 2$)$] If $\grade(I,M)>0$ then 
$$r(I^{n+1},M)=I\cdot r(I^n,M)+\ann M\;\;\mbox{for all}\;n\gg 0.$$
\end{itemize}

Next we relate $\widetilde{I^nM}$ and $I^nM$. In the case when $M = R$, the  following result is proved  in 
\cite[2.13(a)]{Sch}. 
\begin{proposition}\label{impcor1}
Let $I$ be an ideal of $R$ and $M$ an $R$-module. Then
$$\widetilde{I^nM}=I^nM+H^0_I(M)\quad\mbox{for all}\;n\gg 0.$$
\end{proposition}
\begin{proof}  By Proposition \ref{recursion1}(d), $\;\grade(J,N)>0$. Therefore  
$\widetilde{J^nN} = J^nN$\; for all\; $n\gg 0.$ So
\begin{equation*}
\widetilde{J^nN} = \frac{I^nM+H^0_I(M)}{H^0_I(M)}\quad \mbox{for all}\; n\gg 0.\tag{i}
\end{equation*}
It is easy to see that 
$H^0_I(M)\subseteq \widetilde{I^nM}\quad\mbox{for all}\; n\ge 1.$
By an argument similar to Proposition \ref{recursion1}(e), we get 
\begin{eqnarray*}
\widetilde{J^nN} &=& \frac{\widetilde{I^nM}}{H^0_I(M)}.
\end{eqnarray*}
Thus from equation (i), the result follows. 
\end{proof}

\begin{corollary}\label{neoth7}
Assume that ${\widetilde{\mathcal R}}(I,M)=\bigoplus_{n\ge 0} {\widetilde {I^nM}}$ is a Noetherian $\mathcal R(I)$-module and $M$ is separated with respect to the $I$-adic topology. Then $H^0_I(M)=0$.
\end{corollary}
\begin{proof} Since ${\widetilde{\mathcal R}}(I,M)$ is Noetherian, there exists a positive integer $n_0\in \mathbb N$ such that $\widetilde {I^nM}=I^{n-n_0}\widetilde {I^{n_0}M}$ for all $n\geqslant n_0$. So $H^0_I(M)\subseteq I^{n-n_0}\widetilde {I^{n_0}M} \subseteq I^{n-n_0}M$ for all  $n\geqslant n_0$. Thus the result follows from our hypothesis on $M$.
\end{proof}

\section{Relation with a Superficial Element }
\noindent In this section we assume that $(R,\m)$ is local with the maximal ideal $\m$. The goal of this section is to understand the relation between $r(I,M)$ and a superficial element. To ensure the existence of superficial element we assume (unless stated otherwise) that the residue field  $K=R/{\m}$ is infinite. When $I$ is $\m$-primary  our techniques yield an easily computable bound  on $k$ such  that 
$\widetilde{I^n}= (I^{n+k}: I^k)$\; for all $n\ge 1$. 

\s Recall an element $x\in I$ is called $M$-${\it superficial}$ with respect to $I$ if there exists an integer
$c\ge 0$ such that
$$(I^{n+1}M:_M x)\cap I^cM=I^nM\quad \mbox{for all}\; n\ge c.$$
Superficial element exists when $K$ is infinite. If $\grade{(I,M)}> 0$ then every $M$-superficial element is also $M$-regular. Also if $x\in I$ is $M$-superficial and $M$-regular then 
$$(I^{n+1}M:_M x)=I^nM\quad \mbox{for all}\; n\gg 0.\quad\mbox{(see \cite[p. 8]{SJ} for the case $M=R$)}$$

\begin{proposition}\label{8}
Let $x\in I$ be a $M$-superficial element. Then 
\begin{eqnarray*}
 \left( r(I^{n+1},M): x\right)  &=& r(I^n,M) \quad \text{for all } n\ge 1.
\end{eqnarray*}
\end{proposition}
\begin{proof}
 Since $x\in I$ is $M$-superficial, there exists $c>0$ such that 
\begin{equation*}
(I^{j+1}M:_M x)\cap I^cM=I^jM \quad \text{for all } j\ge c.\tag{i}
\end{equation*}
It is easy to see that {$r(I^n,M)\subseteq \left( r(I^{n+1},M): x\right) $}\quad for all $n\ge 1$. Conversely let $a\in \left( r(I^{n+1},M): x\right) $. Then
\begin{eqnarray*}
ax &\in & r(I^{n+1},M)= (I^{n+k+1}M:I^kM)\quad \text{for}\;k\gg 0\\
\mbox{Thus}\;\;axI^kM &\subseteq & I^{n+k+1}M\quad \text{for}\;k\gg 0.
\end{eqnarray*}
Now for $k\gg 0$, we have $aI^kM \subseteq  I^cM$. Therefore, by (i), we get
\begin{eqnarray*}
aI^kM &\subseteq & I^{n+k}M.
\end{eqnarray*}
Thus \; $a \in  (I^{n+k}M:I^kM)\subseteq r(I^n,M)$.
\end{proof}

\s If $\grade{(I, M)}>0$ then $\widetilde{I^nM}=I^nM$ for all $n\gg 0$. Define 
$$\rho^I(M):=\min\{\;n:\;\widetilde{I^iM}=I^iM \;\mbox{for all}\;i\ge n\;\}.$$
If $x\in I$ is $M$-superficial then it is proved in \cite[Corollary 5.3]{TJP} that
$$\rho^I(M)=\min\{\;n:\;(I^{i+1}M:_M x)=I^iM \;\mbox{for all}\;i\ge n\;\}.$$

\begin{proposition}\label{24}
Let $\grade(I,M)>0$. Then 
$$r(I^n, M)=(I^nM: M)\quad\mbox{for all}\; \;n \ge \rho^I(M).$$
\end{proposition}
\noindent{\bf Proof:} Let $x\in I$ be  both $M$-superficial and $M$-regular. So 
\begin{equation*}
(I^{n+1}M:_M x)=I^nM\quad\mbox{for all}\; n\ge \rho^I(M).\tag{i}
\end{equation*}
Let $a\in r(I^n, M)$. Then $aI^kM \subseteq  I^{n+k}M \;\; \mbox{for some}\; k\ge 1$. So $ax^kM \subseteq  I^{n+k}M.$
By repeated use of (i), we get $aM \subseteq  (I^{n+k}M :_M x^k)= I^nM$. Therefore $a\in (I^nM: M)$.
\hfill {$\Box$}

\begin{corollary}\label{22}
Let $x\in I$ be $M$-superficial such that $x^*$ is $G_I(M)$-regular. Then
$$r(I^n, M)=(I^nM: M)\quad \mbox{for all}\;\;n\ge 1.$$
\end{corollary}
\begin{proof}
Since $x^*$ is $G_I(M)$-regular, we have
$$(I^{n+1}M:_M x)=I^nM\quad \mbox{for all}\;n\ge 1$$
So $\rho^I(M)=1$. Therefore the result follows from Proposition \ref{24}. 
\end{proof}

\noindent{\bf Notation:} Let  $I$ be an $\m$-primary ideal. Suppose $M$ is of dimension $d\geq 0$. We define the postulation number $\eta^I(M)$ of $I$ with respect to $M$ as follows:
$$\eta^I(M):=\min\{n\;:\;H_I^M(t)=p_I^M(t),\;\mbox{for all}\;t\ge n\},$$
where $H_I^M(n)=\ell(M/{I^{n+1}M})$ is the Hilbert-Samuel function of $M$ with respect to $I$ and $p_I^M(t)$ is the Hilbert-Samuel polynomial of $M$ with respect to $I$. Let $x\in I$ be $M$-superficial. We set 
$$\eta^I(x,M)=\;\max\{\;\eta^I(M),\;\eta^{I}({M/{xM}})\;\},$$

The following proposition is proved by J. Elias in \cite[1.3]{EJ} for $M=R$. The same proof applies to the general case.  
\begin{proposition}\label{19}
Let $I$ be an $\m$-primary ideal. Let $x\in I$ be an $M$-superficial and $M$-regular element. Then
$$(I^{k+1}M: x)=I^kM\quad\mbox{for all }\;k\geqslant \eta^I(x,M)+1.$$ \qed
\end{proposition}

\begin{remark}\label{obs}
 From above Proposition \ref{19}, if $\grade(I,M)>0$, it follows that $\rho^I(M)\leqslant \eta^I(x,M)+1$. 
\end{remark}

\begin{proposition}\label{23first}
Let $\grade{(I, M)}>0$.  Then 
$$r(I^n, M)=(I^{n+k}M: I^kM)\quad \mbox{for all }\;\;n\geqslant 1\;\;\mbox{and} \;\;k\geqslant \;\rho^I(M).$$
\end{proposition}
\begin{proof}
 Clearly $\grade(I,I^kM)>0$ for all $k$. By Proposition \ref{3}(g) and Proposition \ref{relationf}, we have for each $n\ge 1$,
\begin{eqnarray*}
r(I^n,M) &=& r(I^n,I^kM)\quad \mbox{for all}\;\; k\ge 1,\\
&=& (\widetilde{I^{n+k}M}:I^kM)\quad \mbox{for all}\;\; k\ge 1.\\
\mbox{Thus}\;\;\;r(I^n,M) &=& ({I^{n+k}M}:I^kM)\quad \mbox{for all}\;\; k\ge \rho^I(M).
\end{eqnarray*}
\end{proof}

When $M =R$ we obtain that

\begin{corollary}\label{firstresult}
 Let $I$ be a regular ideal. Then for each value of $n$, we have
$$\widetilde{I^n}=(I^{n+k}: I^k)\;\; \mbox{for all } \;k\geqslant\rho^I(R).$$
\end{corollary}
\begin{remark}
In particular, if $I$ is an $\m$-primary regular ideal then for each $n\ge 1$,
$$\widetilde{I^n}=(I^{n+k}: I^k)\;\; \mbox{for all } \;k\geqslant  \eta^I(x,R)+1.$$
For $n= 1$ the result above was proved by Elias \cite[p.\ 722]{EJ}. However our result does not follow from it.
Furthermore even for $n =1$ our method is simpler to compute.
\end{remark}

It is of interest to find a similar bound for $\widetilde{I^nM}$. We prove

\begin{theorem}\label{23}
Let $\grade{(I, M)}>0$.  Then for each $n\ge 1$, we have
$$\widetilde{I^nM}=(I^{n+k}M:_M I^k)\;\; \mbox{for all }\; k\geqslant \;\rho^I(M).$$
\end{theorem}
\begin{proof}
 Let $x\in I$ be $M$-superficial. It is enough to show that for $k\ge \rho^I(M)$, we have
$$ (I^{n+k+1}M:_M I^{k+1})\subseteq (I^{n+k}M:_M I^k).$$
Let $m\in (I^{n+k+1}M:_M I^{k+1})$. Then $\;mI^{k+1} \subseteq  I^{n+k+1}M$, so \;$mxI^k \subseteq  I^{n+k+1}M$. Therefore
\begin{eqnarray*}
mI^k &\subseteq & (I^{n+k+1}M:_M x)=I^{n+k}M.
\end{eqnarray*}
So $\;m\in (I^{n+k}M:_M I^k)$.
\end{proof}

The next theorem deals with the situation when residue field $R/{\m}$ is not necessarily infinite.

\begin{theorem}\label{21} 
Let $(R,\m)$ be local, $M$ an $R$-module and let $I$ be an ideal with $\grade(I,M)>0$. Then the following hold for all $k\ge\;\rho^I(M),$
\begin{itemize}
\item [$($\rm a$)$] $r(I^n, M)=(I^{n+k}M: I^kM)\;\; \mbox{for all}\;\; n\ge \;1.$
\item [$($\rm b$)$] $\widetilde{I^nM}=(I^{n+k}M:_M I^k)\;\; \mbox{for all }\;n\ge 1.$
\end{itemize}
In particular, when $M=R$ we have for each n,  $\widetilde{I^n}=(I^{n+k}: I^k)\;\; \mbox{for all}\;\;k\ge\;\rho^I(R).$
\end{theorem}
\begin{proof}
 Note that (b) follows from (a). We now prove part (a). Consider the faithfully flat extension 
$$\;R\longrightarrow T=R[X]_{\m R[X]}.$$  
Note that the residue field of $T\;$  is $K(X)$, the quotient field of polynomial ring $K[X]$ and it is infinite. Set $q=IT$ and $E=M\otimes_R T$. By Proposition \ref{3}(d), we get
\begin{eqnarray*}
r(I^n,M)\otimes_R T &=& r(q^n, E)\\
&=& (q^{n+k}E:_T q^kE)\quad \mbox{for all}\;k\ge \rho^q(E),\\
&=&(I^{n+k}M:I^kM)\otimes_R T\quad \mbox{for all}\;k\ge \rho^q(E).
\end{eqnarray*}
By \cite[1.7]{TJP}, we have  $\;\rho^q(E)=\rho^I(M)$. Fix $k\ge \rho^I(M)$ and set $$D=r(I^n,M)/{(I^{n+k}M:I^kM)}.$$ Then we have $\;D\otimes_R T=0.$ Since $T$ is a faithfully flat extension of $R$, we get $D=0$.
\end{proof}

We used the packages {\it CoCoA} \cite{cocoa}\; and {\it Singular} \cite{Sin}\; for our computations.
We reconsider the  Example \ref{equality}. In this example we apply Theorem \ref{21} to compute $\widetilde{I^n}$ for each $n$. 

\begin{example}\label{exampleone}
(see \cite[1.4]{RS}) Let $R=k[x,y]_{\left\langle x,y \right\rangle }$ and the ideal
\begin{eqnarray*}
I &=& \; \left\langle y^{22},x^4y^{18},x^7y^{15},x^8y^{14},x^{11}y^{11},x^{14}y^8,x^{15}y^7,x^{18}y^4,x^{22}\right\rangle R.
\end{eqnarray*}
Set $u=x^{22}+y^{22}$. The Poincare series of $I$  and $I/{u}$ are
\begin{eqnarray*}
PS_I(t) &=& \frac{227+189t+10t^2+10t^3-2t^4}{{(1-t)}^2},\\
PS_{I/{u}}(t) &=& \frac{227+189t+12t^2+6t^3}{{(1-t)}}.
\end{eqnarray*}
So $e_i(I/{u})=e_i(I)$\; for $i=0,1$. Thus $u$ is $R$-superficial with respect to $I$. Note that $\rho^I(R)\le \eta^I(x,R)=3$. Therefore by \ref{21}, $\widetilde {I^n} = (I^{n+3}:I^3)=I^n\;\;\mbox{for all} \; n\geqslant 3.$ Also 
\begin{eqnarray*}
\widetilde {I} &=& (I^4:I^3)=I,\\
\widetilde {I^2} &=&  (I^5:I^3)=I^2+ \left\langle x^{20}y^{24},x^{24}y^{20} \right\rangle R.
\end{eqnarray*}
\end{example}

\section{The case when $I$ has a Principal Reduction }
\noindent In this section we discuss  ideals having principal reductions. When $I$ has a principal reduction $J=(x)$, the computation of $r(I,M)$ is greatly simplified. Let $r=r_x(I)=\min\{\;n:\;I^{n+1}=xI^n\;\}$, the reduction number of $I$ with respect to $(x)$ then we show that $r(I^n,M)=(I^{n+r}M:I^rM)$\; for all $n\geqslant 1$. This result is then used to compute many examples.

\begin{proposition}\label{rednum}
Let $\grade(I,M)>0$. If $I$ has a principal reduction $(x)$ with reduction number $r_x(I) \leq r$ then the following hold
\begin{itemize}
\item [$($\rm a$)$]  If $\;r(I^{n+1}, M)=x\cdot r(I^n, M)$ for some n then $r(I^{m+1}, M)=x\cdot r(I^m, M)$ for all $\;m\ge n$. 
\item [$($\rm b$)$] $r(I^n, M)= (I^{n+r}M:I^rM)$\; for all\; $n\ge 1$.
\item [$($\rm c$)$] $\widetilde{I^nM}=(I^{n+r}M:_M I^r)$ \; for all \;$n\ge 1$.
\end{itemize}
\end{proposition}
But before we do this we first need to prove the following lemma.
\begin{lemma}\label{15}
Let the situation be as in Proposition \ref{rednum} . Then $x$ is $M$-regular and
$$(x)\cap r(I^{n+1}, M)=x\cdot r(I^n, M)\;\; \mbox{for all}\;\;n\ge 1.$$
\end{lemma}
\begin{proof}
 It is easy to check that $x$ is $M$-regular. Let $\;ax\in r(I^{n+1}, M)$ for some $a\in R$. Then 
$axI^kM \subseteq  I^{n+1+k}M\quad\mbox{for some}\;k\ge 0.$
We assume $k\ge r_x(I)$. Therefore
\begin{eqnarray*}
axI^kM &\subseteq & xI^{n+k}M.
\end{eqnarray*} 
Since $x$  is $M$-regular, we get  $a\in (I^{n+k}M:I^kM)\subseteq r(I^n, M).$
\end{proof}

We now give the proof of proposition.

\begin{proof}[Proof of Proposition \ref{rednum}]
 {\bf (a)} By Lemma \ref{15}, $x$ is $M$-regular. For all $m\ge n$, we have
\begin{eqnarray*}
r(I^{m+1}, M) &\subseteq & r(I^{n+1}, M)\subseteq (x).\\
\mbox{Thus,}\;\;\;\;\;r(I^{m+1}, M)&\subseteq &  (x)\cap r(I^{m+1}, M), \\
&=& x\cdot r(I^m, M), \;\;\;\;\;\;\mbox{(by Lemma \ref{15})}\\
&\subseteq & r(I^{m+1}, M).
\end{eqnarray*}
(b): It is enough to show that $(I^{n+r+1}M:I^{r+1}M) \subseteq (I^{n+r}M:I^rM)$. Let $a\in (I^{n+r+1}M:I^{r+1}M)$.
\begin{eqnarray*}
aI^{r+1}M &\subseteq & I^{n+r+1}M,\\
\mbox{so}\quad ax\cdot I^rM &\subseteq & x \cdot I^{n+r}M.
\end{eqnarray*} 
As $x$  is $M$-regular we get $aI^rM \subseteq   I^{n+r}M.$\\
\noindent {\bf (c)} It is similar to (b).
\end{proof}

For $M = R$ Proposition \ref{rednum}(c) yields
\begin{corollary}
Let $I$ be a regular ideal having a principal reduction $(x)$ with  reduction number $r_x(I) \leq r$. Then
$\widetilde{I^n}=(I^{n+r} \colon I^r)$ \; for all \;$n\ge 1$. \qed
\end{corollary}

\begin{remark}
For $n = 1$ the result above has been proved already in \cite[2.1]{HGD}.
However for $n \geq 2$ their result does not imply ours. 
 Notice that, since $r_{x^n}(I^n) \leq r_x(I) \leq r$, \cite[2.1]{HGD}
implies $\widetilde{I^n} = (I^{nr + n} \colon I^{nr})$.
\end{remark}

\noindent{\bf Examples:}
We use Proposition \ref{rednum} to construct many examples. We first give an example of a module $M$ with no free summand such that $r({I}^n,M)={I}^n$ \; for all\; $n\geqslant 1\;$ but $\;{IM}\ne \widetilde{IM}.$

\begin{example}\label{singular1}
Let $\;Q=k[x,y]_{\left\langle x,y\right\rangle }\;$ be a local ring with the maximal ideal $\n$. Set $(R,\m)=\left( Q/{\left\langle y^3\right\rangle },\n/{\left\langle y^3\right\rangle }\right) $. 
Consider the $R$-module 
$$M= \left\langle \left(\begin{array}{cc} 0 \\  y^2  \end{array} \right),\left(\begin{array}{cc} y \\  x  \end{array} \right) \right\rangle \;\subseteq\; R^2 .$$ 
Note that $M$ is  a Cohen-Macaulay $R$-module. Also notice that $\;{\m}^3=x\cdot {\m}^2\;$ so $\;r_x(\m)=2$.  We  compute $r({\m}^n,M)$ for $n\ge 1$. Using Proposition \ref{rednum}(b), one  checks 
\begin{eqnarray*}
r({\m}^i,M)=({\m}^{i+2}M:{\m}^2M) &=& {\m}^i\quad \mbox{for }\;i=1,2,3.
\end{eqnarray*}
Thus, by Proposition \ref{rednum}(a) we have 
$$ r({\m}^{n+1},M) = x\cdot r({\m}^n,M)={\m}^{n+1}\;\mbox{for all}\; n\ge 2.$$ 
We also compute $\widetilde{\m M}$ by using Proposition \ref{rednum}(c)
\begin{eqnarray*}
\widetilde{\m M} &=& ({\m}^3M:_M {\m}^2) = \left\langle \left  (\begin{array}{cc} 0 \\ y^2 \end{array}\right), \left  (\begin{array}{cc} y^2 \\ xy \end{array}\right),\left  (\begin{array}{cc} xy \\ x^2 \end{array}\right)\right\rangle \\
\mbox{and}\quad {\m}M &=& \left\langle \left (\begin{array}{cc}   y^2 \\ xy \end{array} \right), \left (\begin{array}{cc}   xy \\ x^2 \end{array} \right)\right\rangle. 
\end{eqnarray*}
Hence\; ${\m}M \ne  \widetilde{\m M}.$
\end{example}

We give an example where
$r({I}^n,M)=\widetilde{{I}^n}\;$ for all \;$n\ge 1$ and ${I M}= \widetilde{I M}.$ 

\begin{example}
Let $$A=\frac{k[x,y,z]}{<x^5-z^2,y^3-xz>}\simeq k[t^6.t^7,t^{15}].$$ 
Set $R=A_{\n}$, where $\n =\left\langle x,y,z\right\rangle A$. Let $\m$ be the maximal ideal of $R$.
Let $$M=\left\langle \left  (\begin{array}{cc} x^2 \\ z^2 \end{array} \right),\left  (\begin{array}{cc} y \\ yz \end{array} \right) \right\rangle \;\subseteq \; R^2.$$
Note that $M$ is a Cohen-Macaulay $R$-module and $\;{\m}^6=x\cdot {\m}^5\;$ so $\;r_x(\m)=5$. Therefore using Proposition \ref{rednum}(b), one can check the following
\begin{eqnarray*}
r({\m}^n,M)=({\m}^{n+5}M:{\m}^5 M) &=& \widetilde{\m^n} =({\m}^{n+5}:{\m}^5) \;\; \mbox{for}\;\; n=1,...,6.\\
\mbox{Also}\quad \widetilde{\m^n} \ne  {\m}^n \quad \mbox{for}\;\; n=2,3,4
&\mbox{and}& \widetilde{\m^n} =  {\m}^n\;\; \mbox{for all}\; n\ge 5.
\end{eqnarray*}
Thus, by Proposition \ref{rednum}(a) we have 
$$ r({\m}^{n+1},M) = x\cdot r({\m}^n,M)={\m}^{n+1}\;\;\mbox{for all}\;\; n\ge 5.$$ 
Using Proposition \ref{rednum}(c), one can verify that $\widetilde{\m M}=({\m}^6M:_M{\m}^5)=\m M$.
\end{example} 
 
Next we give an example in which we have $\widetilde{I^{n}}\ne r({I}^n,M)$\; for all\; $n\geqslant 1$, and ${\mathcal F}^I_M$ is a stable $I$-filtration.
 
\begin{example}\label{examplerepeat1}
Let $A = k[t^4,t^{11}, t^{17},t^{18}]$. Using Singular we get
\begin{eqnarray*} 
A &\simeq &B=\frac{k[x,y,z,w]}{\left\langle y^2-xw,yz-x^7,z^2-x^4w,yw-x^3z,zw-x^6y,w^2-x^2yz \right\rangle }.
\end{eqnarray*} 
Set $R=B_{\m}$, where $\m =\left\langle x,y,z,w\right\rangle B$.
Let \;$I=\left\langle x,y \right\rangle R$\;\;  and\; $M=\left\langle x,y,z \right\rangle R$.
One can check $\;r_x(I)=2$. Note that $x$ is both $M$-regular and $R$-regular. By Proposition \ref{rednum}(b), we have 
\begin{eqnarray*}
r(I,M) &=& (I^3M:I^2M)=\left\langle x,y,z,w  \right\rangle R,\\
r(I^2,M) &=& (I^4M:I^2M)=\left\langle x^2,xy,xz,xw  \right\rangle R=x\cdot r(I,M).\\
\mbox{And}\quad\widetilde I &=&(I^3:I^2)=\left\langle x,y,w\right\rangle R,\\
 \widetilde {I^2}&=&(I^4:I^2)=\left\langle x^2,xy,xw\right\rangle R=x\cdot \widetilde I.
\end{eqnarray*} 
Therefore we have
\begin{eqnarray*}
r(I^{n+1},M) &=& x\cdot r({I}^n,M)=\left\langle x^{n+1},x^ny,x^nz,x^nw \right\rangle R\;\;\mbox{for all}\;\; n\ge 1,\\
\widetilde{I^{n+1}}=&=& x\cdot \widetilde{I^{n}}=\left\langle x^{n+1},x^ny,x^nw\right\rangle R\;\;\mbox{for all}\;\; n\ge 1.
\end{eqnarray*}
Notice that $x^{n-1}z\leftrightarrow t^{4n+13}$.  We use the expression of $I^n$ from Example \ref{9} to verify that\; $x^{n-1}z\notin \widetilde{I^{n}}.$ 
Therefore  $\;\widetilde{I^{n}}\ne r({I}^n,M)$\; for all\; $n\geqslant 1$. Since $M$ is a non-zero ideal of $R$ and $R$ is a domain so $\ann M=0$. Notice that $I$ is a regular ideal. So, by Theorem \ref{stablity2},  ${\mathcal F}^I_M$ is a stable $I$-filtration.
\end{example}

\section{Application to Hilbert Functions}
\noindent In this section we assume that $(R,\m)$ is a local ring with the maximal ideal $\m$ and $I$ an $\m$-primary ideal. Let $P^M_I(t)$ be the Hilbert-Samuel function of $M$ with respect $I$. We first show that if $\grade(I,M)>0$ then the set
$$\mathcal{H}(I)=\{\;J: J\;\mbox{is an ideal of}\;R\;\mbox{such that}\;J\supseteq I\;\mbox{and}\;P^M_J(t)=P^M_I(t)\;\}$$ 
has $r(I,M)$ as the  unique maximal element. If $\dim M=1$ and $\depth M=0$ then ${\widetilde{\mathcal R}}(I,M)=\bigoplus_{n\ge 0} {\widetilde {I^nM}}$ is not Noetherian. However $\widetilde G_I(M)=\bigoplus_{n\geqslant 0} \widetilde{I^n}/{\widetilde{I^{n+1}}}$ is a Noetherian $G_I(R)$-module, in fact a Cohen-Macaulay $G_I(R)$-module (Theorem \ref{neoth7}). Next we give an application of Proposition \ref{impcor1} to show that if  $\dim M=1$ then 
$e^I_1(M)-e^I_0(M)+\ell(M/{IM})\geqslant -\ell\left( H^0_{\m}(M)\right).$ 

\s Recall that the Hilbert-Samuel function of $M$ with respect $I$ is the function 
$$n\longmapsto \ell\left( M/{I^{n+1}M}\right) \quad \text{for all}\  n\gg 0.$$
It is well known that for all large values of $n$ it is given by a polynomial $P^M_I(n)$ of degree $r=\dim M$,  the {\it Hilbert-Samuel polynomial} with respect to $I$. It can be written in the form
\[
 P^{M}_{I}(X) = \sum_{i=0}^{r}(-1)^i e_{i}^{I}(M)\binom{X + r - i}{r-i}.
\]
The integers $e^{I}_{0}(M), e^{I}_{1}(M),...,e^{I}_{r}(M)$ are called the {\it Hilbert coefficients} of $M$ with respect to $I$. The number $e^{I}_{0}(M)$ is also called the {\it multiplicity} of $M$ with respect to $I$.
 
\s  For an $\m$-primary ideal $I$, we define the set
$$\mathcal{H}(I):=\{ J \colon J\ \ \text{is an ideal of}\ R \ \mbox{such that}\ J\supseteq I  \ \ \mbox{and}\ P^M_J(t)=P^M_I(t) \ \}.$$

\begin{proposition}\label{uniqueratliff}
Let $\grade(I,M)>0$. Then $r(I,M)$ is the unique maximal element of  $\;\mathcal{H}(I)$.
\end{proposition}
\begin{proof} By Corollary \ref{25}, we have ${r(I,M)}^nM=I^nM\;\;\mbox{for all}\;\;n\gg 0.$
So $r(I,M)\in \mathcal{H}(I)$. 
Conversely if $J\in \mathcal{H}(I)$ then $J^nM=I^nM$\;\;for all\; $n\gg 0$.
Let $x\in J$. We have 
\begin{eqnarray*}
xI^{n-1}M = xJ^{n-1}M &\subseteq& J^nM=I^nM\;\;\mbox{for all}\;n\gg 0.\\
\mbox{So}\quad x &\in & (I^nM: I^{n-1}M)=r(I,M).
\end{eqnarray*}
Thus $J\subseteq r(I,M)$. It follows that $r(I,M)$ is the unique maximal element in $\mathcal{H}(I)$.
\end{proof}

\begin{remarks}\label{impremarks}
We have 
\begin{itemize}
\item [$($\rm a$)$] $\widetilde{I^nM}=I^nM+H^0_{\m}(M)\quad \mbox{for all}\;\;n\gg 0.$ (by Proposition \ref{impcor1})
\item [$($\rm b$)$] By using Artin-Rees lemma, we can check that
$$I^nM\cap H^0_I(M)=0\quad \mbox{for all}\;\;n\gg 0.$$
\item [$($\rm c$)$] Using (a) and (b), it is easy to prove that
$$\frac{\widetilde{I^nM}}{I^nM} \simeq \frac{H^0_I(M)}{I^nM\cap H^0_I(M)} \simeq   H^0_I(M)\quad \mbox{for all}\;\;n\gg 0.$$
\end{itemize}
\end{remarks}

Set ${\widetilde G_I}(M)=\bigoplus_{n\ge 0}{\widetilde {I^nM}}/{\widetilde {I^{n+1}M}}$. In dimension one we have

\begin{theorem}\label{finiteness2}
Let $(R,\m)$ be a local ring. Let $I$ be an $\m$-primary ideal and suppose that $M$ is an $R$-module with $\dim M=1$. Then
\begin{itemize}
\item [$($\rm a$)$]  $\ell(M/{\widetilde {I^{n+1}M}})= e^I_0(M)(n+1)-e^I_1(M)- \ell \left(H^0_{\m}(M)\right)$.
\item [$($\rm b$)$]  ${\widetilde G_I}(M)$ is a finitely generated $ G_I(R)$-module of dimension 1.
\item [$($\rm c$)$]  ${\widetilde G_I}(M)$ is a Cohen-Macaulay $ G_I(R)$-module.
\end{itemize}
\end{theorem}

\begin{proof} We may assume that $K=R/{\m}$ is infinite, otherwise we consider the standard base change $R\longrightarrow R[X]_{\m A[X]}$ (see \cite[1.3]{TJP}).\\
\noindent  {\bf (a)} From the exact sequence 
$$0 \longrightarrow {\widetilde{I^{n+1}M}}/{I^{n+1}M} \longrightarrow M/{I^{n+1}M}\longrightarrow M/{\widetilde {I^{n+1}M}} \longrightarrow 0$$
and using Remark \ref{impremarks}(c), it follows that for all $n\gg0,$
\begin{eqnarray*}
\ell(M/{\widetilde {I^{n+1}M}}) &=& e^I_0(M)(n+1)-e^I_1(M)- \ell \left(H^0_{\m}(M)\right).
\end{eqnarray*}
\noindent{\bf (b)} Consider the exact sequence 
$$0 \longrightarrow {\widetilde{I^{n}M}}/{\widetilde {I^{n+1}M}} \longrightarrow M/{\widetilde {I^{n+1}M}}\longrightarrow M/{\widetilde {I^{n}M}} \longrightarrow 0.$$
By (a), we have  
\begin{equation*}
\ell\left( {\widetilde{I^nM}}/{\widetilde {I^{n+1}M}}\right) = e^I_0(M),\quad\mbox{for all}\;\;n\gg 0.\tag{ii}
\end{equation*}
Let $x\in I$ be $M$-superficial which exists as the residue field $K$ is infinite. Let $x^*$ be the image of $x$ in $I/{I^2}$, considered as $G_I(R)$-element. Then by Proposition \ref{8}, the sequence 
\begin{equation*}
0 \longrightarrow {\widetilde{I^{n}M}}/{\widetilde {I^{n+1}M}} \xrightarrow{x*}  {\widetilde{I^{n+1}M}}/{\widetilde {I^{n+2}M}}\tag{iii} 
\end{equation*}
is exact. In view of (ii), we get
\begin{equation*}
{\widetilde{I^{n+1}M}}/{\widetilde {I^{n+2}M}}\simeq  x^* \cdot \left( {\widetilde{I^{n}M}}/{\widetilde {I^{n+1}M}}\right) \quad\mbox{for all}\;\;n\gg0, \;\mbox{say for}\;n\ge n_0.\tag{iv} 
\end{equation*}
 Let $S=\{\;z_1,z_2,...,z_s\;\}$ be the set of homogeneous generators of $\bigoplus_{i=0}^{n_0} {{\widetilde{I^{n}M}}}/{{\widetilde {I^{n+1}M}}}$ as an $R$-module. Using (iv) it follows that $S$ generates ${\widetilde G_I}(M)$ as a $G_I(R)$-module. Since $\ell\left( {\widetilde{I^nM}}/{\widetilde {I^{n+1}M}}\right)=e_0^I(M)>0$ for all $n\gg 0$. It follows that $\dim {\widetilde G_I}(M)=1.$\\
\noindent {\bf (c)} By(iii), we have $x^*$ is ${\widetilde G_I}(M)$-regular. So $\depth {\widetilde G_I}(M)>0$. As $\dim {\widetilde G_I}(M)=1$ we have ${\widetilde G_I}(M)$ is Cohen-Macaulay.
\end{proof}

\begin{theorem}\label{inequality1}
Let the situation be as in Theorem \ref{finiteness2}. Then 
$$e^I_1(M)-e^I_0(M)+\ell(M/{IM})\geqslant -\ell \left( H^0_{\m}(M)\right).$$
Further equality holds if and only if $\;\widetilde{IM}=IM$ and $\widetilde{I^nM}=I^nM+H^0_{\m}(M)$ for $n\geqslant 2$.
\end{theorem}
In \cite[3.1]{GN} a sharper bound is found. However the case of equality is different.
\begin{proof} As ${\widetilde G_I}(M)$ is a finitely generated $G_I(R)$-module so  there exists $h(t)\in \mathbb Z[t]$ such that
$$\sum_{n=0}^{\infty} \ell\left( {\widetilde{I^{n}M}}/{\widetilde {I^{n+1}M}}\right) t^n =\frac{h(t)}{(1-t)},$$
where $h(t)=h_0+h_1t+h_2t^2+...+h_st^s.$ As ${\widetilde G_I}(M)$ is a Cohen-Macaulay $G_I(R)$-module, all the coefficients of $h(t)$ are non-negative. Thus
\begin{eqnarray*}
{\widetilde {e_1}^I}(M) &= & {\widetilde{e_0}}^I(M)-h_0+\sum_{j= 2}^s (j-1)h_j,\\
&=& {\widetilde{e_0}}^I(M)- \ell\left(M/{\widetilde{IM}}\right)+\sum_{j=2}^s (j-1)h_j,\\
&= &  {\widetilde{e_0}}^I(M)-\ell\left( M/{IM}\right)+\ell \left( {\widetilde{IM}}/{IM}\right) +\sum_{j=2}^s (j-1)h_j.
\end{eqnarray*}
For $n \gg 0$ we have $\ell\left( M/{\widetilde{I^{n+1}M}}\right) = {\widetilde {e_0}}^I(M)(n+1)-{\widetilde {e_1}}^I(M)$. By comparing it with  \ref{finiteness2}(a), we get ${\widetilde{e_0}}^I(M) = e_0^I(M) \;\;\mbox{and}\;\;{\widetilde {e_1}}^I(M) = e_1^I(M)+\ell\left( H^0_{\m}(M)\right).$ Thus we get
\begin{eqnarray*}
e_1^I(M) &=& {\widetilde {e_1}}^I(M)-\ell\left( H^0_{\m}(M)\right), \\
&= &{\widetilde{e_0}}^I(M)-\ell\left( M/{IM}\right) +\ell \left( {\widetilde{IM}}/{IM}\right)+\sum_{j=2}^s (j-1)h_j -\ell\left( H^0_{\m}(M)\right), \\
&=& e_0^I(M)-\ell\left( M/{IM}\right)+\ell \left( {\widetilde{IM}}/{IM}\right)+\sum_{j=2}^s (j-1)h_j  -\ell\left( H^0_{\m}(M)\right).
\end{eqnarray*}
Hence equality holds if and only if $\;\sum_{j=2}^s (j-1)h_j=0$\; and \;${\widetilde{IM}}=IM$.
Assume that equality holds that is, ${\widetilde{IM}}=IM$ and $\;\sum_{j=2}^s (j-1)h_j=0$. We have $h_j=0$ for $j\ge 2$. So, by graded Nakayama's Lemma one can see that  ${\widetilde G_I}(M)$ is generated in degree $0$ and $1$. This follows that $\widetilde{I^jM}=I^jM+\widetilde{I^{j+1}M}$\; for all $j\geqslant 2$. By Proposition \ref{impcor1} we thus get $\widetilde{I^jM}=I^jM+H^0_{\m}(M)$ for $j\ge 2$. Conversely suppose that $\widetilde{IM}=IM$ and $\widetilde{I^jM}=I^jM+H^0_{\m}(M)$ for $j\geqslant 2$. This gives $\widetilde{I^jM}=I^jM+\widetilde{I^{j+1}M}$\; for all $j\geqslant 2$. Which implies that ${\widetilde G_I}(M)$ is generated in degree zero and one. So we have $h_j=0$ for all $j\ge 2$.
\end{proof}

We give an example which shows that the bound in Theorem \ref{inequality1} can be attained.

\begin{example}\label{examfirst10}
Let $R=k[|x,y|]/(xy,y^2)$ with the maximal ideal $\m= \;<x,y>R$. The Hilbert function of $R$ with respect to $\m$ is 
$$\frac{1+z-z^2}{1-z}.$$
We get $e_1(R)-e_0(R)+\ell(R/{\m})=-1$. We now claim $\ell(H^0_{\m}(R))=1$. Let $q=<y>R$ be the ideal in $R$. Notice that $\m\cdot q=0$. Also $q\ne 0$. Thus  $q\simeq k$ as $R$-modules. Further $R/q\simeq k[|x|]$ is Cohen-Macaulay. Now using the short exact sequence 
$$0\rightarrow q\rightarrow R \rightarrow R/q \rightarrow 0,$$
we get $H^0_{\m}(R)\simeq H^0_{\m}(q)=k$. 
\end{example}

\indent In view of  Theorem \ref{inequality1}\; we pose a question\\
\noindent{\bf Question:} If $\dim M\geqslant 2\;$ then {\it is the set
$$\{\;e^I_1(M)-e^I_0(M)+\ell(M/{IM})\;:\; I\;\;an \;\m\mbox{-primary}\; ideal\;\}$$
bounded below ?}\\
This result holds if $M$ is a generalized Cohen-Macaulay $R$-module. (see \cite[2.4 and 3.1]{GN}).

\section*{Acknowledgments}
It is a pleasure to thank Prof. J. K. Verma and Prof A. V. Jayanthan for many discussions.
 We also thank the referee for many pertinent comments.


\providecommand{\bysame}{\leavevmode\hbox to3em{\hrulefill}\thinspace}
\providecommand{\MR}{\relax\ifhmode\unskip\space\fi MR }
\providecommand{\MRhref}[2]{%
  \href{http://www.ams.org/mathscinet-getitem?mr=#1}{#2}
}
\providecommand{\href}[2]{#2}

\end{document}